\documentclass[12pt]{article}

\usepackage{flushend, cuted}
\usepackage{amsmath}
\usepackage{amsfonts}
\usepackage{bbm}

\usepackage{amssymb}
\usepackage{float}
\usepackage{color,amsmath,amssymb, amsfonts,amstext,amsthm}

\usepackage{epsfig, graphicx, graphics}
\usepackage{longtable}
\usepackage{cite}
\usepackage{subfigure}

\usepackage{amsfonts}
\usepackage{amssymb}
\usepackage{color, amsmath,amssymb, amsfonts, amstext,latexsym}

\usepackage{amssymb, epsfig, amssymb, latexsym}
\usepackage{amsmath}
\usepackage{graphicx}
\usepackage{longtable}

\newcommand{\s}{\sigma}

\renewcommand{\phi}{\varphi}

\newcommand{\R}{{\mathbb R}}

\newtheorem{example}{Example}

\title{{\bf Detecting Stochastic Governing Laws with Observation on Stationary Distributions }\footnote{This work is supported by the National Natural Science Foundation of China (NSFC) (Grant No.11901536). Xiaoli Chen is supported by the Ministry of Education, Singapore, under its Research Centre of Excellence award to the Institute for Functional Intelligent Materials (I-FIM, project No. EDUNC-33-18-279-V12)}}
\author{\centerline{\bf Xiaoli Chen$^{a,b,}\footnote{xlchen@nus.edu.sg}$,
Hui Wang$^{c,}\footnote{Corresponding author: huiwang2018@zzu.edu.cn}$,
Jinqiao Duan$^{d,}\footnote{duan@iit.edu}$}
\centerline{${}^a$Department of Mathematics}
\centerline{National University of Singapore, 119077, Singapore}
\centerline{${}^b$Institute for Functional Intelligent Materials}
\centerline{National University of Singapore, 117544, Singapore}
\centerline{${}^c$ School of Mathematics and Statistics} \centerline{Zhengzhou University, Zhengzhou, China}
\centerline{${}^d$ Department of Applied Mathematics \& Department of Physics,} \centerline{Illinois Institute of Technology, Chicago, IL 60616, USA}}

\begin{document}

\maketitle

\begin{abstract}
Mathematical models for complex systems are often accompanied with uncertainties. The goal of this paper is to extract a stochastic differential equation governing model with observation on stationary probability distributions.
We develop a neural network method to learn the drift and diffusion terms of the stochastic differential equation. We introduce a new loss function containing the Hellinger distance between the observation data and the learned stationary probability density function. We discover that the learnt stochastic differential equation provides a fair approximation of the data-driven dynamical system after minimizing this loss function during the training method. The effectiveness of our method is demonstrated in numerical experiments.

\end{abstract}




\section{Introduction}
Mathematical models for scientific and engineering systems often involve uncertainties and thus are often in the form of stochastic differential equations (SDE). These stochastic dynamical systems are ubiquitous in biology, physics, geosciences and other fields. Stochastic dynamical systems provide an appropriate framework to investigate random phenomena \cite{Arnold2003, Pas, Sura,GaoDuan, Wu}. Hence, the determination of SDE models is crucial for quantifying and predicting dynamical behaviors of the nonlinear system under random fluctuations.
An SDE is characterized by the drift term and diffusion terms. In this paper, we aim to detect appropriate drift and diffusion terms with stationary probability distribution data.

The stationary probability distribution of a stochastic dynamical system does not change with time and it is the stationary solution of the corresponding Fokker-Planck equation \cite{Hung, Hairapetian, Yarmchuk, Gefen, ArnoldWishtutz, Liberzon, Gray}.
The stationary probability distribution carries the information of the underlying stochastic system \cite{Khasminskii, Schmalfuss}.
The Hellinger distance is the distance between probability distributions which characterizes how close two different distributions are.
In this paper, we use Hellinger distance to identify whether the constructed SDE is an appropriate
approximation of a data-driven stochastic dynamical system.


Neural networks can be represented as compositions of simple functions with parameters, and such functional representations can be used for parameter estimation of time-series data and kernel estimation \cite{Gerber}.
There has been some progress in learning stochastic differential equation models from noisy data.
A variation estimation method  was used to learn the drift term  with the observation trajectory data \cite{Batz, batz2018approximate, opper2017estimator, opper2019variational}.
 There was also an RNN-based variational method \cite{ryder2018black}, a sparse learning method \cite{boninsegna2018sparse}, and a Kramers-Moyal formulae \cite{tabar2019analysis} for learning stochastic dynamical systems.
A stochastic adjoint sensitivity method was proposed to learn stochastic differential equations \cite{li2020scalable} or stochastic differential equations  with jumps \cite{jia2019neural}.
In \cite{YangLiu-GAN-SODE}, they used small samples from just a few snapshots of unpaired data to infer the drift and diffusion terms of stochastic differential equations. Moreover, in \cite{Zhang2020, Xu2020}, they learned  L\'evy noise parameters by deep neural networks. In \cite{Zhang2022}, they solved the steady-state Fokker-Planck equation with a small amount of data through combining the deep KD-tree.

 \medskip

We have recently developed a data-driven approach \cite{li2021data,Lu} to discover stochastic  differential equations
with non-Gaussian L\'{e}vy noise using the nonlocal Kramers-Moyal formulas, and further
  learned the stochastic differential equations from discrete particle samples
at different time snapshots using the Fokker-Planck equation and physics-informed neural networks \cite{chen2021solving}.


However, in addition to sample path observation data, there are recent advances in observing or measuring stationary probability distributions \cite{ Hairapetian, Gefen, Liberzon,Yang_I}. To take advantage of these new types of data, we devise a neural network method to extract stochastic dynamical system models with stationary probability distribution or a long time trajectory as observation data. This motivates our research reported in this paper. Specifically, we develop a neural network method to extract stochastic governing laws based on probability measures. Given observation data, we learn the drift and diffusion terms which are approximated by two neural networks. Since if we learn the drift and diffusion together, the results would not be unique. So in this work we proposed two approaches. The first approach entails simply learning the drift or diffusion terms. The second technique involves learning the drift and diffusion terms simultaneously with one drift term observational data. We compare our learned results in three-dimensional settings with Hellinger distance substituted by Jensen-Shannon divergence and mean-square distance which demonstrate the efficacy of our proposed approaches.

This paper is organized as follows. In Section \ref{sec:2}, we present two methods for learning stochastic governing laws based on physics informed neural networks and Hellinger distance of probability distributions. In Section \ref{sec:3}, we present examples to learn the drift terms and the diffusion terms. Finally, we end with some discussions in Section \ref{sec:4}.

\section{Methodology}\label{sec:2}
\label{methods}

\subsection{Problem setup }
Consider the following stochastic differential equation (SDE)
\begin{equation}\label{sde3456}
    dX_t=b(X_t)dt + \sigma(X_t) dB_t, \;\; X_0=x_0,
\end{equation}
where the n-dimensional vector function $b(\cdot)$ is the drift term, the $n\times n$ matrix function $\sigma(\cdot)$ is the diffusion term, and $B_t $ is an n-dimensional Brownian motion.

The generator of the SDE \eqref{sde3456} is \cite{Duan2015}:
\begin{align}  \label{generator}
  Au=&\sum_{i=1}^{n} b_i\frac{\partial u}{\partial x_i}+\frac{1}{2}\sum_{i,j=1}^{n}(\sigma \sigma^T)_{i,j}\frac{\partial^2 u}{\partial x_i \partial x_j}. \nonumber
 \end{align}

The probability density function (PDF) is a quantity that carries information of the stochastic system. The  time evolving probability density function of the solution process $X_t$ is governed by the Fokker-Planck equation, which is written as follows:
\begin{equation}
    \begin{aligned}
&\partial_t p(x,t)=A^{*}p(x,t), ~~~~~x \in \mathbb{R}^n, t>0,  \\
&p(x,0) = p_0(x),
 \end{aligned}
 \end{equation}
where $p_0(x)$ is the initial probability density function, $A^{*}$ is the adjoint operator of the generator $A$ and has the following form:
\begin{equation}\label{ajoint}
    \begin{aligned}
  A^{*}p=&-\sum_{i=1}^{n}\frac{\partial }{\partial x_i}(b_i p)+\frac{1}{2}\sum_{i,j=1}^{n}\frac{\partial^2}{\partial x_i \partial x_j}((\sigma \sigma^T)_{i,j} p).
 \end{aligned}
 \end{equation}
Note that the Fokker-Planck equation is a deterministic linear partial differential equation with an initial condition.

The stationary Fokker-Planck equation is:
\begin{equation}\label{stationary-FPE}
    \begin{aligned}
  A^{*}p(x)=&0,
 \end{aligned}
 \end{equation}
 with a condition  $ \int_{\R^n} p(x)dx=1$. We assume that there exists a unique stationary probability density function (still denote it by $p(x)$) in this paper.

We consider the scenario when available data is  time series observation data or probability density function. Our objective is to infer the drift and diffusion terms. Because of the stochasticity  of the dynamical system, we could not use the mean square error to get the loss function of the SDE \eqref{sde3456}. The main issue is how to  quantify the stochastic dynamics with use the deterministic indexes. For example, we can use the maximal likelihood estimation \cite{felix_SDE} or the most transition pathway \cite{chen_learnMTP} to extract or learn the SDE model. Here we will use the stationary probability density function as the deterministic index to learn  the SDE. If the available data is  long time trajectory data of $X(t)$, we may first use kernel density estimation to learn the probability density function.
We will propose a machine learning method to learn the drift and the diffusion terms of the SDE, with
different measures for the distance of the observed probability density function.


\subsection{Machine Learning}

As the drift $b$ and diffusion $\s$ characterize the uncertainty of the SDE,  we will  estimate them based on observations of probability distributions (i.e., probability measures) of the system paths $X_t$.
Now we  introduce the Hellinger distance \cite{Cha,Beran} between two probability distributions. It is  used to quantify the distance between two probability distributions in the space of probability measures. For our purpose here, the Hellinger distance $H(p, q)$ between two probability density functions $p(x) $ and $q(x)$ is defined as follows,
\begin{equation}\label{Hellinger}
   H(p, q) \triangleq  \sqrt{\frac12 \int (\sqrt{p(x)} - \sqrt{q(x)} \; \;)^2 dx},
\end{equation}
which satisfies the property: $0 \leq H(p, q) \leq 1$.


With the observed stationary probability density $q(x)$, we determine or estimate the drift term $b(x)$  by minimizing the Hellinger distance between the true stationary probability density $p(x)$ for the solution process $X(t)$ and the observation probability density  $q(x)$.




Note that  Hellinger distance is a measure to describe the distance of two probability density functions. Other distance also can describe the distance. Such as,
given the probability density function $p(x)$ and $q(x)$, respectively, the Kullback-Leibler (KL) divergence is defined as
\begin{equation}\label{eqn:KL divergence}
    \begin{aligned}
    H_{KL}(p||q) =\int p(x)\log(\frac{p(x)}{q(x)}) dx.
    \end{aligned}
\end{equation}
While the Kullback-Leibler divergence is asymmetry, there also exists a symmetric  measure between two probability density, which is Jensen-Shannon divergence, introduced as follows:
\begin{equation}\label{eqn:JS divergence}
    \begin{aligned}
    H_{JS}(p||q) =\frac{1}{2}H_{KL}(p||\frac{q+p}{2})+\frac{1}{2}H_{KL}(\frac{q+p}{2}||p).
    \end{aligned}
\end{equation}
Later  we will also use the Jensen-Shannon divergence to measure the distance.


Given the noise intensity $\sigma(x)$ and observation of the stationary probability density $q(x)$, we will learn the drift term $b$. We devise two neural networks to approximate the drift term and stationary probability density $p(x)$, where the input is the space domain $x$ and the output is the $b_{NN}(x)$ and $p_{NN}(x)$.

On the one hand, the output of $p_{NN}(x)$ should satisfy the functional \eqref{Hellinger2}. We define the loss function as:
\begin{equation} \label{loss_OM}
Loss_H=\frac{1}{2N_H}\sum_{i=1}^{N_H}(\sqrt{p_{NN}(x_i)}-\sqrt{q(x_i)})^2,
\end{equation}
where $\{x_i\}_{i=1}^{N_H}$ are the points in the spatial domain to compute the integral and $N_H$ is the number of the observation data.

On the other hand, the neural networks of the drift term $b_{NN}(x)$ and the stationary probability density $p_{NN}(x)$ should satisfy the steady Fokker-Planck equation \eqref{stationary-FPE}.

\begin{figure}[h]
\centerline{\includegraphics[width=8cm]{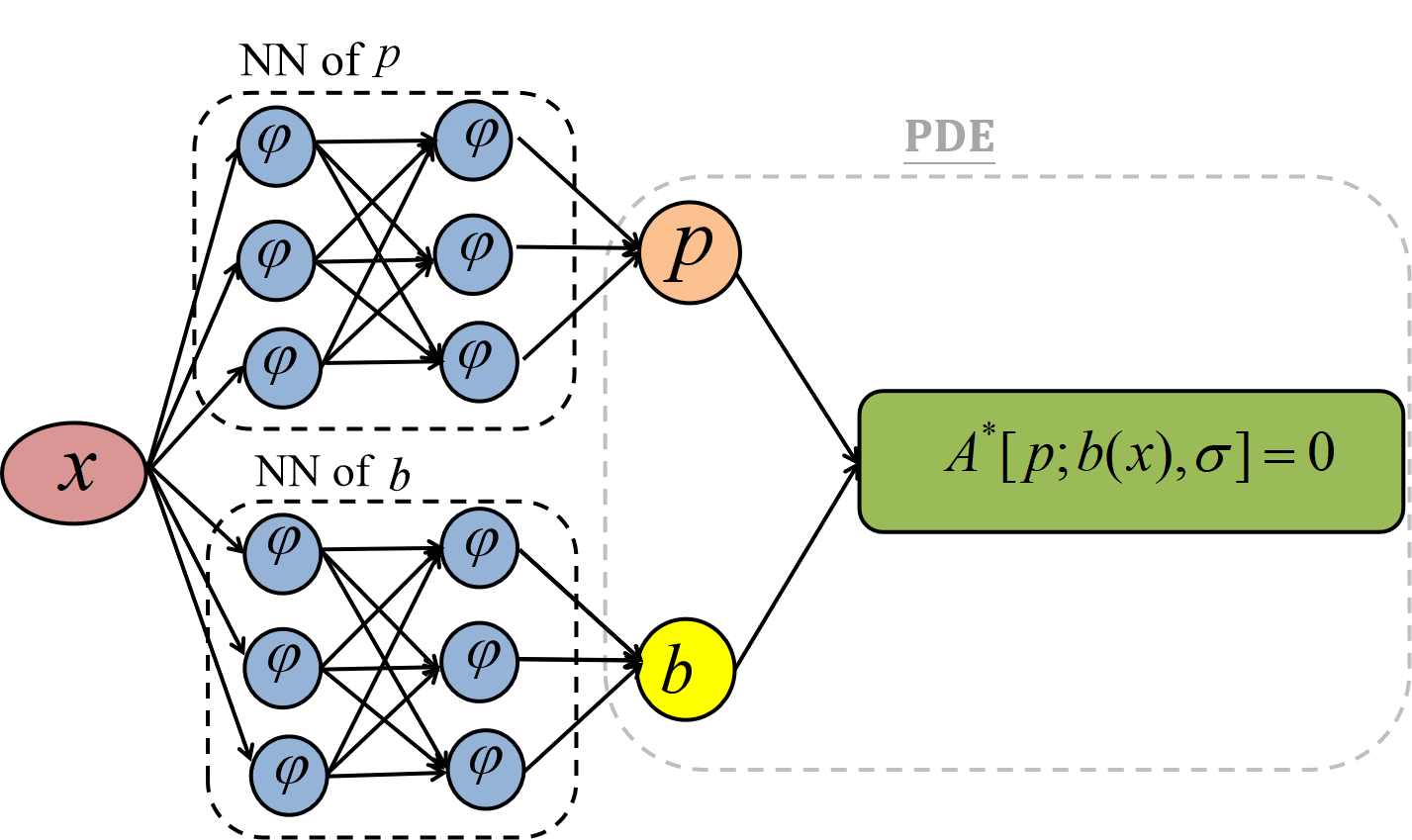}}
\caption{\textbf{Schematic of the neural network for solving PDEs}:
two neural networks for approximate the probability p and the drift term b, where the input is x and $\phi$ is the activation function.}
\label{NN_sketch}
\end{figure}

Similar to the physics informed neural network \cite{Raissi2019,chen2020}, we define the residual neural network as
\begin{equation}\label{1d_rhs}
f(x)=-\sum_{i=1}^{n}\frac{\partial }{\partial x_i}(b_{iNN} p_{NN})+\frac{1}{2}\sum_{i,j=1}^{n}\frac{\partial^2}{\partial x_i \partial x_j}((\sigma \sigma^T)_{i,j} p_{NN}).
\end{equation}
Then the loss function of the residual neural network is defined as:
\begin{equation} \label{loss_f}
Loss_{f}=\frac{1}{N_f}\sum_{i=1}^{N_f}(f(x_i))^2,
\end{equation}
where $\{x_i\}_{i=1}^{N_f}$ is the residual points in the spatial domain and $N_f$ is the number of the residual points. Here we randomly choose the residual points at each iteration step.  The sketch of the method is shown in Figure \ref{NN_sketch}.

The total loss function   is
\begin{equation} \label{loss}
Loss=Loss_{H}+Loss_{f}.
\end{equation}

 For the unknown drift term and diffusion term, because $b(x)=0$ and $\sigma(x)=0$ is also the minimization solution of loss function, thus we could not learn the terms uniquely if we train the loss function \eqref{loss}.
The observation data \cite{chen2020} of drift term at some points need to know. To avoid the zeros solution, the observation data of drift term at few points is given, i.e. $\{x_i,b(x_i)\}_{i=1}^{N_b}$. The loss function of the drift term is:
\begin{equation} \label{loss_drift}
Loss_{b}=\frac{1}{N_b}\sum_{i=1}^{N_b}(b_{NN}(x_i)-b(x_i))^2,
\end{equation}
so the loss function is defined as:
\begin{equation} \label{loss_case2}
Loss_2=Loss_{H}+Loss_{f}+Loss_b.
\end{equation}

If we use the Jensen-Shannon divergence to measure the distance, we just replace the $loss_H$ to
\begin{equation} \label{loss_JS}
Loss_{JS}=H_{JS}(p_{NN}||q).
\end{equation}

We also compare our method \eqref{loss} with the traditional physics informed neural network (PINN) method \cite{Raissi2019} for solving the inverse problem of the Fokker-Planck equation. With the observation data of the stationary probability density function $q(x)$, the loss function is written as mean square error:
\begin{align}\label{loss_ob}
Loss_{ob}=\frac{1}{N_{ob}}\sum_{i=1}^{N_{ob}}(p_{NN}(x_i)-q(x_i))^2.
\end{align}
The total loss of PINN method is defined as:
\begin{align}\label{loss_PINN}
Loss_{PINN}=Loss_{ob}+Loss_{f}.
\end{align}

If the observed data is the trajectory of stochastic differential equation, i.e., $X=(X_{t_0}, X_{t_1},\cdots, X_{t_N})$, we can use the kernel density estimation to obtain the probability density function, denoting as $q_{KD}$. Similarly, we  replace the probability density function loss $loss_H$ in Eq. \ref{loss_OM} or $loss_{JS}$ in Eq. \ref{loss_JS} with another loss from the estimated density, as $q_{KD}$.


Remark: For an SDE with non-Gaussian  L\'evy case
\begin{align}  \label{stomodel_levy}
  dX_t &= b(X_t)dt +  \varepsilon dL_{t}^{\alpha},  ~~X_t \in \mathbb{R}^n,
 \end{align}
where $b(\cdot)$ is the vector drift term, and $L_t^\alpha$ is a symmetric $\alpha$-stable L\'evy process in $\mathbb{R}^n$. The generating triplet of the  L\'evy process is $(0,0,\nu_\alpha)$.\\
The corresponding   nonlocal Fokker-Planck operator is \cite{Duan2015}
\begin{equation}\label{ajoint-levy}
    \begin{aligned}
  A^{*}p=&-\sum_{i=1}^{n}\frac{\partial }{\partial x_i}(a_i p)+\varepsilon^\alpha \int_{\mathbb{R}^n\setminus \{0\}}[p(x+y)-p(x)]\nu_{\alpha}(dy),
 \end{aligned}
 \end{equation}
 where $\nu_{\alpha}(dy)$ is the $\alpha$-stable L\'evy measure and $\nu_{\alpha}(dy)=C_{n,\alpha}||y||^{-n-\alpha}dy$, $C_{n,\alpha}=\frac{\alpha \Gamma((n+\alpha)/2) }{2^{1-\alpha} \pi^{n/2}\Gamma(1-\alpha/2)}.$  The stationary probability density function
 is the solution for the nonlocal equation $A^{*}p=0$.
To use our method to learn the SDE \eqref{stomodel_levy} driven by L\'evy noise, we only need to change the loss function of residual neural network \eqref{loss_f} to \eqref{ajoint-levy}. As for the nonlocal integral term, we   discretize it with a scheme in our earlier work \cite{chen2019most}, and while for  the first order derivative term, we evaluate with automatic differentiation \cite{auto}.

\section{Numerical Experiments}\label{sec:3}

We first present an analytical example to learn a simple stochastic system, with quite involved calculations. For more complex stochastic systems, we will have to use our proposed machine learning method as demonstrated in the following numerical experiments.

\subsection{Analytical method for learning stochastic dynamical systems}

Consider a scalar stochastic differential equation
\begin{equation}
    dX_t = b(X_t) dt +\sigma dB_t,
\label{sde3456b}
\end{equation}
with appropriate conditions on drift $b$ and diffusion $\s$ (see \cite[p.170]{Klebaner}), such as, $b \leq 0$ and $\s \neq 0$ as well as some smoothness requirements, there exists a unique stationary probability density $p(x)$ for the SDE \eqref{sde3456b}, as a solution of the steady Fokker-Planck equation,
\begin{equation} \label{p}
 p(x)= \frac{C}{\s^2(x)} e^{\int_{x^*}^x \frac{2b(y)}{\s^2(y)}dy},
 \end{equation}
where the positive normalization constant $C$ is chosen so that $p > 0$ and $ \int_{\R} p(x)dx=1$, i.e.,
\[C \triangleq 1/ \int_{-\infty}^{\infty} \frac{e^{\int_{x^*}^x \frac{2b(y)}{\s^2(y)}dy}}{\s^2(x)} \; dx.   \]
Note that $x^*$ here may be an arbitrary reference  point so that the integral $\int_{x^*}^x \frac{2b(y)}{\s^2(y)}dy$ exists. Different choice of $x^*$  only affects the normalization constant $C$. (Say, take $x^*=0$ if that is possible).


Given the observed stationary probability density $q$, we like to find out the true  stationary probability density $p$. Consider   Hellinger distance $H$  between probability  densities $p(x)$ and $q(x)$.
\begin{equation}\label{Hellinger2}
   H(b, \sigma)  \triangleq  \sqrt\frac12 \int_{\R} (\sqrt{p(x)} - \sqrt{q(x)} \; \;)^2 dx,
\end{equation}
where $p(x)$ is in \eqref{p}.
The corresponding Euler-Lagrange equation for $H^2$  is
\begin{equation} \label{Euler}
\frac{d}{dt}H^2_{\sigma}=H^2_b.
\end{equation}
Since the Euler-Lagrange equation is the necessary condition  for the  functional to obtain the  minimum. So we can solve the corresponding Euler-Lagrange equation to get the minimum value of the  functional.

\begin{example}
Consider a specific scalar stochastic model
\[  dX=b(X)dt +  dB_t,  \]
with unknown drift  $b(x)$, and given diffusion $\sigma=1$. Given an ``observation" of the stationary probability density $q(x)= \frac{1}{ \sqrt{2\pi}}  \; e^{- \frac12 x^2 }$ (the Gaussian  distribution).
Find a   function $b(x)$ so that the Hellinger distance $H^2(b(x)) = \frac12 \int_{\R} [\sqrt{p(x)} -\sqrt{q(x)}]^2 dx$ is minimized.
\end{example}

 The true stationary probability density for the solution process $X_t$ is
 \begin{align}\label{SPD}
p(x)=\frac{e^{2\int_{0}^xb(y)dy}}{\int_{-\infty}^{\infty} e^{2\int_{0}^x b(y)dy }dx}.
\end{align}

The Euler-Lagrange  equation is a necessary condition for functional minima
\begin{equation}\label{Hellinger2-1}
   I(b)= \frac12 \int_{\R} (p(x)+q(x)-2\sqrt{p(x)}\sqrt{q(x)}\;) dx.
\end{equation}

Submitting Eq. (\ref{SPD}) and $q(x) $ into Eq. (\ref{Hellinger2-1}), we get
\begin{equation}\label{He}
\begin{split}
   I(b)&= \frac{1}{2\sqrt{2\pi}}\int_{\R}e^{-\frac{1}{2}x^2}dx+\frac{1}{2}\int_{\R}
   \frac{e^{2\int_0^x b(y)dy}}{\int_{\R}e^{2\int_0^x b(y)dy}dx}dx
   -\int_{\R}\sqrt{\frac{e^{- \frac12 x^2+2\int_0^x b(y)dy}}{\sqrt{2\pi}\int_{\R}e^{2\int_0^x b(y)dy}dx}}dx . \\
\end{split}
\end{equation}

In order to get the minima of $I(b)$, we can obtain
$ I^{ \prime}(b)=0$ and $b(x)=-kx, k\geq 0$. \\
Submitting $b(x)$ into $p(x)$, then $p(x)=\sqrt{\frac{k}{\pi}}e^{-kx^2}$, which satisfies $\int_{-\infty}^{\infty}p(x)dx=1$. \\
The error $Err=\|p(x)-q(x)\|_H$, submitting $p(x)$ and $q(x)$ into $Err$:

\begin{equation}
\begin{aligned}
Err=&\|p(x)-q(x)\|_H=\frac{1}{2}\int_R(\sqrt{p(x)}-\sqrt{q(x)})^2dx  \\
=&\frac{1}{4\sqrt{\pi}}+\sqrt{\frac{k}{2\pi}}-\sqrt{\frac{4k}{2\pi(k+\frac{1}{2})}}.\\
\end{aligned}
\end{equation}
This is a function $f(k)$ about variable $k$. $f_{min}$ attains when  $k=\frac{1}{2}$. So, $b(x)=-\frac{1}{2}x.$

In this example, we can luckily find the optimal drift term $b$ analytical. While it is exceedingly difficult to compute the true drift term by hand in many problems. So, in the cases below, we use our proposed machine learning method to learn the drift and diffusion terms.

\subsection{Machine learning for learning stochastic dynamical systems}

The neural networks in our numerical experiments below have 4 hidden layers and 20 neurons per layer, with $\tanh$ activation function. The weights are initialized with truncated normal distributions. The biases are initialized as zero. We use the Adam optimizer with a learning rate $10^{-4}$ to train the loss function.

\begin{figure}[h]
\begin{minipage}[]{0.3 \textwidth}
 \leftline{\tiny\textbf{(a)}}
\centerline{\includegraphics[width=4.5cm]{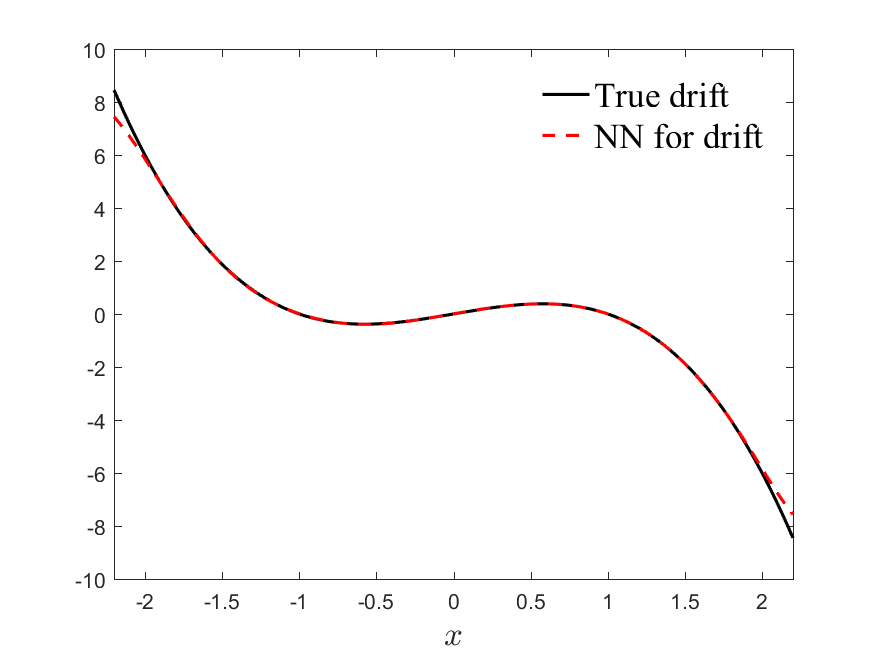}}
\end{minipage}
\hfill
\begin{minipage}[]{0.3 \textwidth}
 \leftline{\tiny\textbf{(b)}}
\centerline{\includegraphics[width=4.5cm]{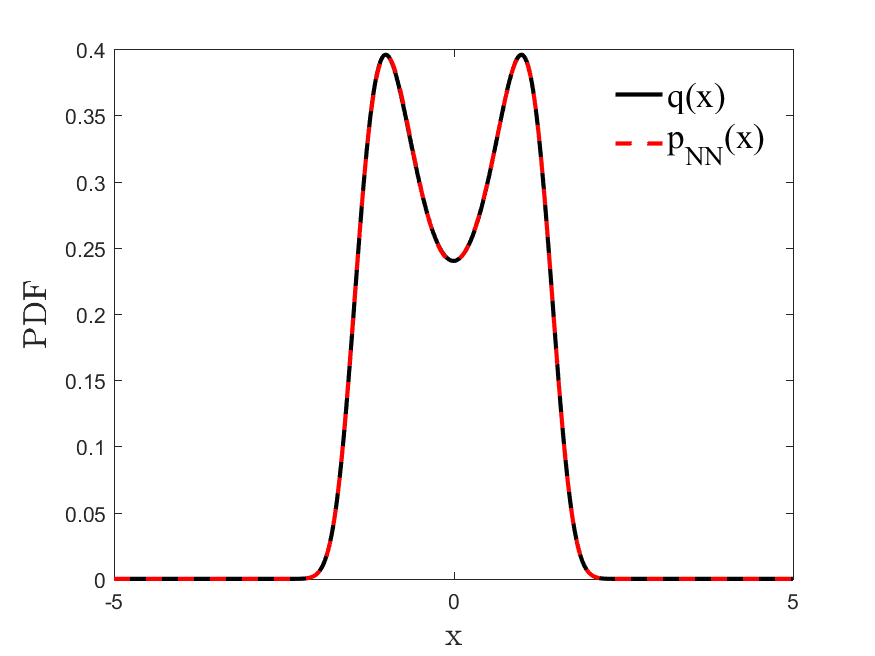}}
\end{minipage}
\hfill
\begin{minipage}[]{0.3 \textwidth}
 \leftline{\tiny\textbf{(c)}}
\centerline{\includegraphics[width=4.5cm]{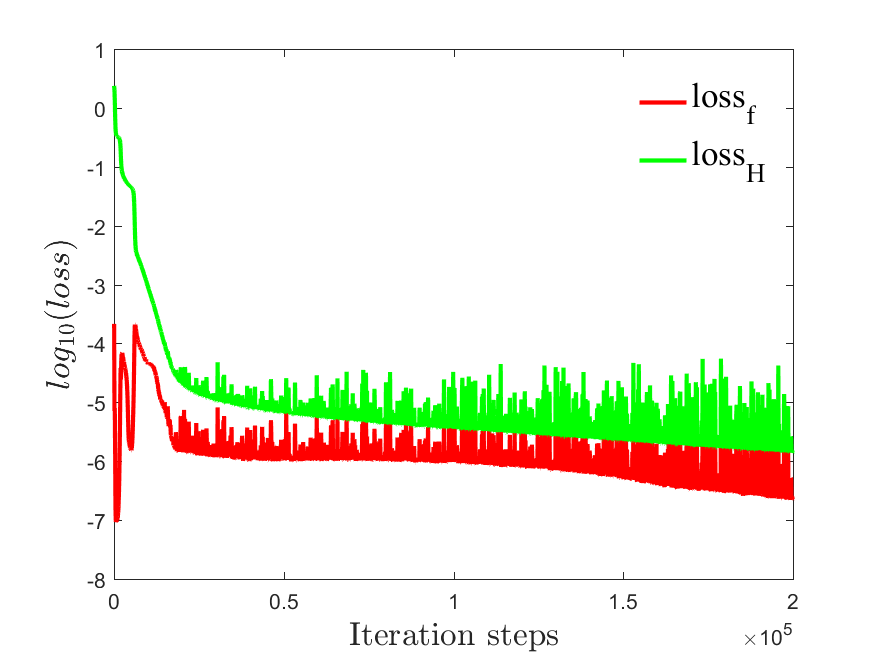}}
\end{minipage}
\caption{\textbf{Unknown drift of Example 2:} (a) learned drift term; (b) learned probability; (c) loss function.}
\label{case2}
\end{figure}

\begin{example}
Consider a scalar stochastic model
\[  dX=b(X)dt + \sigma dB_t,  \]
with drift function $b(x) =x-x^3$. Given an ``observation" of the stationary probability density $q(x)=\frac{1}{A}  e^{\frac{1}{\sigma^2} x^2-\frac{x^4}{2}}$, where $A=\int_{\R}e^{\frac{1}{\sigma^2} x^2-\frac{x^4}{2}} dx$.
Find a drift  function  $b(x)$  so that the Hellinger distance $I(b(x)) = \frac12 \int_{\R} [\sqrt{p(x)} -\sqrt{q(x)}]^2 dx$ is minimized.
\end{example}
Given the noise intensity (diffusion)  $\sigma=1$,   we use two fully connected neural networks to approximate the drift term and stationary probability density respectively. We choose  $N_{H}=1001$, $N_f=10000$ to train the loss function \eqref{loss}. The results we learned are shown in Figure \ref{case2}. In Figure \ref{case2} (a), we plot the true drift term (black line) and the learned drift term (red line). The neural network can approximate the drift very well. In Figure \ref{case2} (b), the given $q(x)$ and neural network result of $p_{NN}(x)$ can  approximate well too. We also plot the loss function evolves with the number of iterative steps. The loss is less than $10^{-4}$.

\begin{figure}[h]
\begin{minipage}[]{0.3 \textwidth}
 \leftline{\tiny\textbf{(a)}}
\centerline{\includegraphics[width=4.5cm]{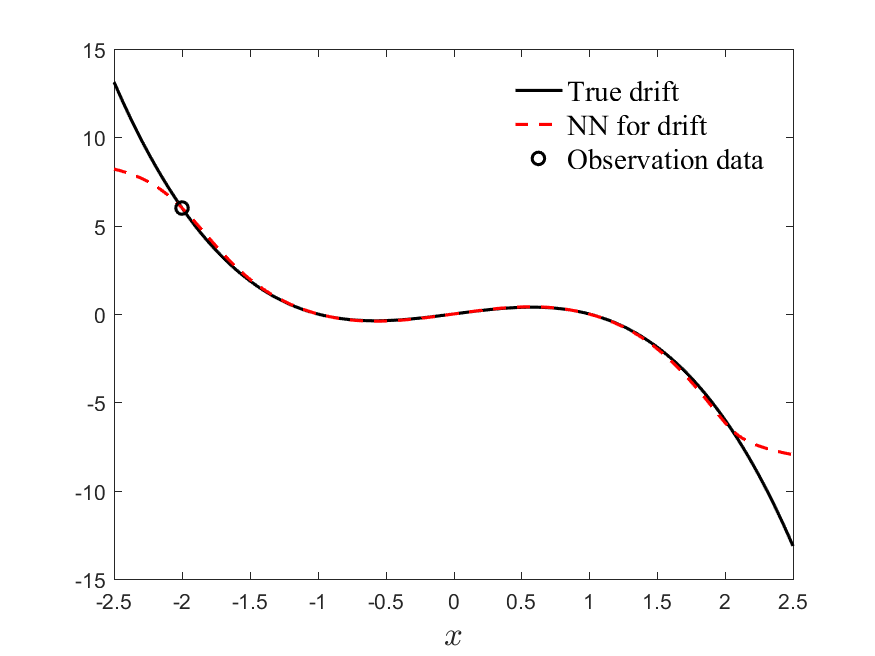}}
\end{minipage}
\hfill
\begin{minipage}[]{0.3 \textwidth}
 \leftline{\tiny\textbf{(b)}}
\centerline{\includegraphics[width=4.5cm]{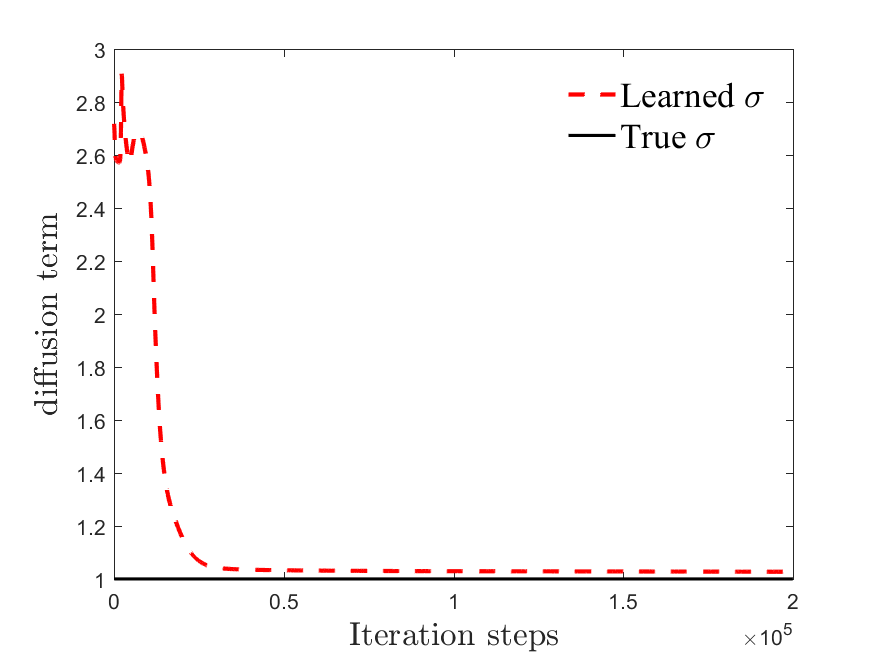}}
\end{minipage}
\hfill
\begin{minipage}[]{0.3 \textwidth}
 \leftline{\tiny\textbf{(c)}}
\centerline{\includegraphics[width=4.5cm]{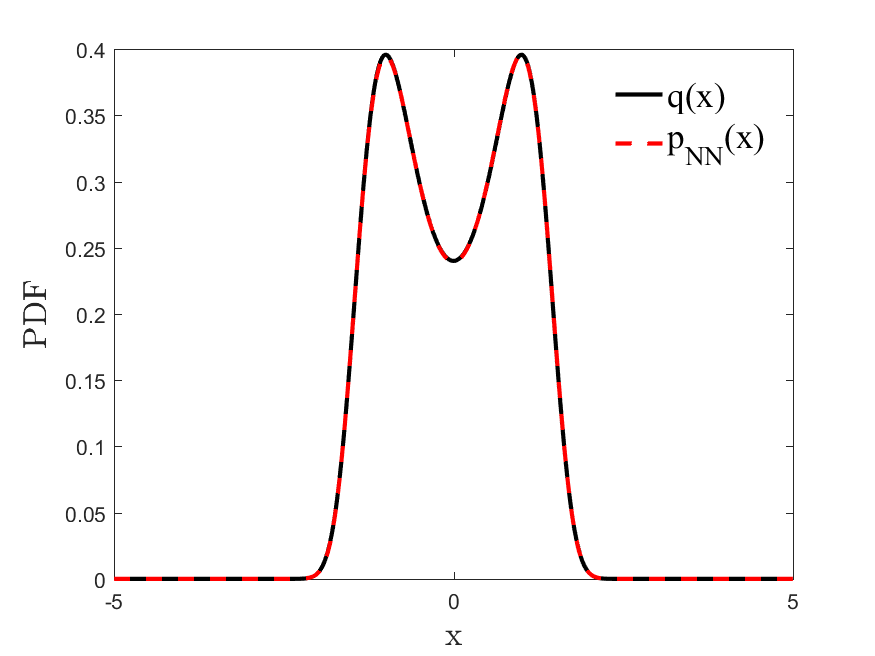}}
\end{minipage}
\caption{\textbf{ Unknown drift and diffusion  of Example 2:} (a) the learned drift term with 1 observation data at $x=-2$; (b) the learned diffusion term; (c)the learned probability density.}
\label{case_dw_diffusion}
\end{figure}

For the unknown drift term and diffusion term, we train the loss function \eqref{loss_case2} to get the optimal results. Only one observation data of the drift term at $x=-2$ is given. The results are shown in Figure \ref{case_dw_diffusion}.
 We present learned drift result in Figure \ref{case_dw_diffusion} (a), and the diffusion term evolution predictions as the iteration of the optimiser progresses in Figure \ref{case_dw_diffusion} (b). The neural network  of probability density is shown in Figure \ref{case_dw_diffusion} (c).
 We can see for one observation data of the drift term, the drift and diffusion term can be learned well. When $|x|>2$, the error of the learned drift term becomes larger than that for the other $x$. The fundamental reason for this is that the probability in bigger $x$ is almost zero, providing very little information there.

\begin{figure}[h]
\centerline{\includegraphics[width=18cm]{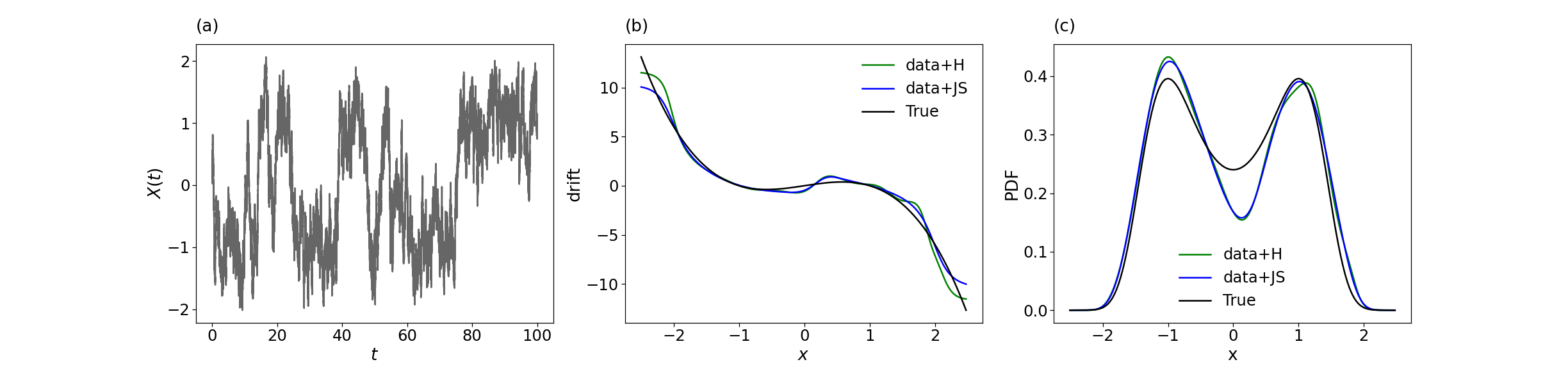}}
\caption{\textbf{ Unknown drift of Example 2 with one trajectory observation data:} (a) trajectory observation data $X(t)$; (b) the learned drift term; (c)the learned probability density.}
\label{case_dw_tra}
\end{figure}

We also use our method to learn the SDE model, with     only   one trajectory observation data $X(t)$. The observation data is the long time trajectory of $X(t)$ and is shown in Figure \ref{case_dw_tra} (a). We   use kernel density estimation to get the probability density and then use our method to learn the drift  and diffusion terms. The results of the drift term are shown in Figure \ref{case_dw_tra} (b) and the probability density is shown in Figure \ref{case_dw_tra} (c). The results validate that our method also works with long time trajectory observation data, while the error of the probability density function is larger than using the stationary probability density function observation data. If we have more trajectory data of the $X(t)$, we can learn the probability density function better.

\begin{example}
Consider a scalar stochastic model
\[  dX=b(X)dt +  \sigma dB_t,  \]
with drift function $b(x)$. Given an ``observation" of the stationary probability density $q(x)=\frac{1}{\pi}\frac{1}{1+x^2}$ .
Find a function $b(x)$ so that the Hellinger distance $I(b(x)) = \frac12 \int_{\R} [\sqrt{p(x)} -\sqrt{q(x)}]^2 dx$ is minimized.
\end{example}

Similar to the last example, we fix the noise intensity (diffusion) $\sigma=1$, and use two neural networks to approximate the drift term and stationary probability density respectively. Here we also choose $N_{H}=1001$, $N_f=10000$.

\begin{equation} \label{p3}
 I(b(x))=\frac12 \int_{\R} (\sqrt{p(x)} -\sqrt{q(x)})^2 dx,
 \end{equation}
where $q(x)=\frac{1}{\pi}\frac{1}{1+x^2}$.

\begin{figure}[h]
\begin{minipage}[]{0.3 \textwidth}
 \leftline{\tiny\textbf{(a)}}
\centerline{\includegraphics[width=4.5cm]{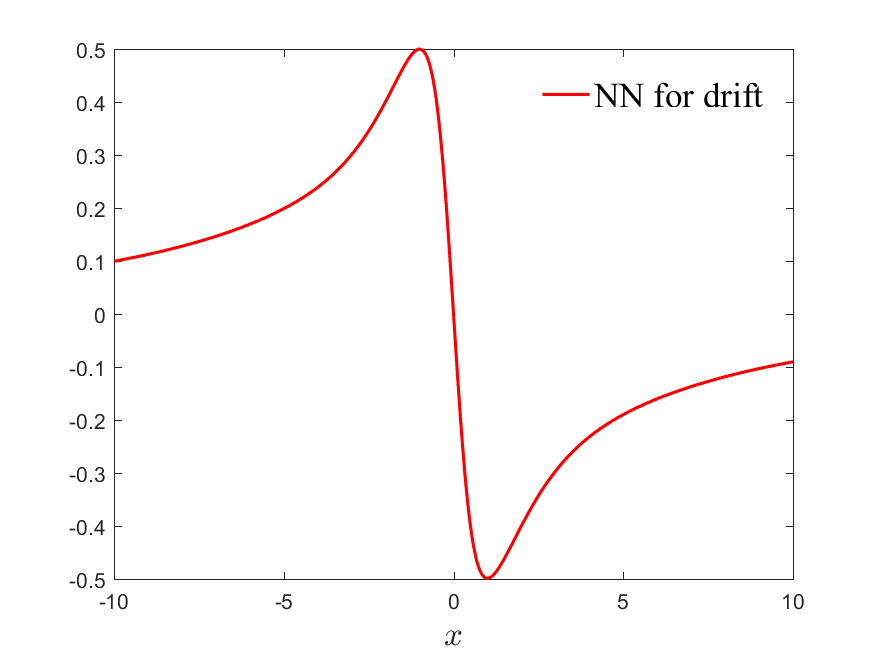}}
\end{minipage}
\hfill
\begin{minipage}[]{0.3 \textwidth}
 \leftline{\tiny\textbf{(b)}}
\centerline{\includegraphics[width=4.5cm]{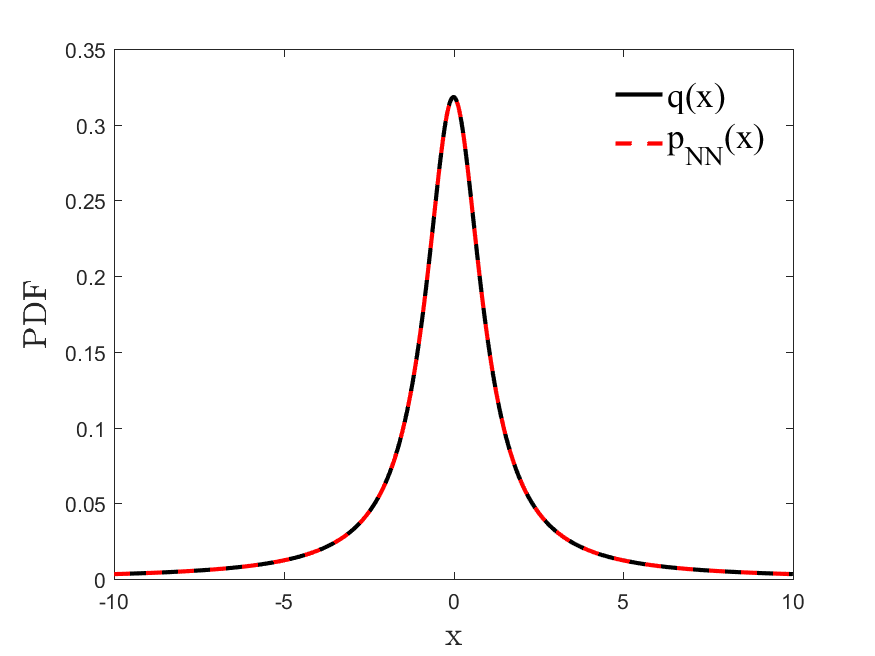}}
\end{minipage}
\hfill
\begin{minipage}[]{0.3 \textwidth}
 \leftline{\tiny\textbf{(c)}}
\centerline{\includegraphics[width=4.5cm]{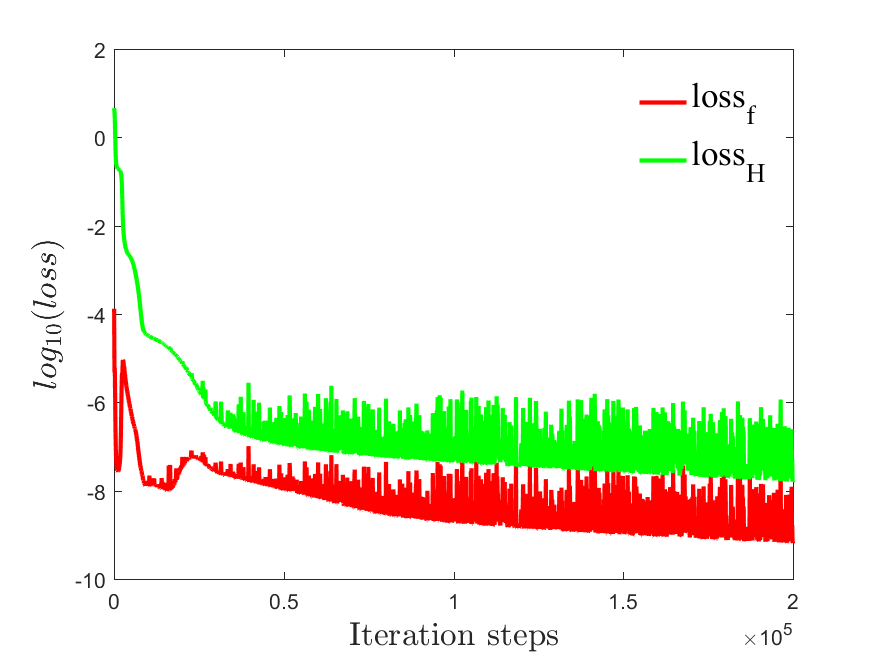}}
\end{minipage}
\caption{\textbf{Unknown drift of Example 3:} (a) learned drift term; (b) learned probability; (c) loss function of OM part $loss_H$ and equation part $loss_f$.}
\label{case3}
\end{figure}

The results are shown in Figure \ref{case3}. In this case, the true drift term $b(x)$ is unknown, so we could not compare the   drift term. By comparing the learned  stationary distribution and the observation data, the result shows that they match well. What is more, the loss function is sufficiently  small, and  the probability density distribution is concentrated around zero. So zero could be the stable point of this system. Our learned drift term $b(x)$  has one stable point zero, as in Figure \ref{case3}. This is in line with our expectations.

Molecular and cell biology is playing an increasingly important role in life sciences. For example,   many research findings are about stochastic fluctuations inducing phenotypic diversity in gene expression. Here we
consider a stochastic gene regulation model \cite{wang1, wang2}.

\begin{figure}[h]
\begin{minipage}[]{0.3 \textwidth}
 \leftline{\tiny\textbf{(a)}}
\centerline{\includegraphics[width=4.5cm]{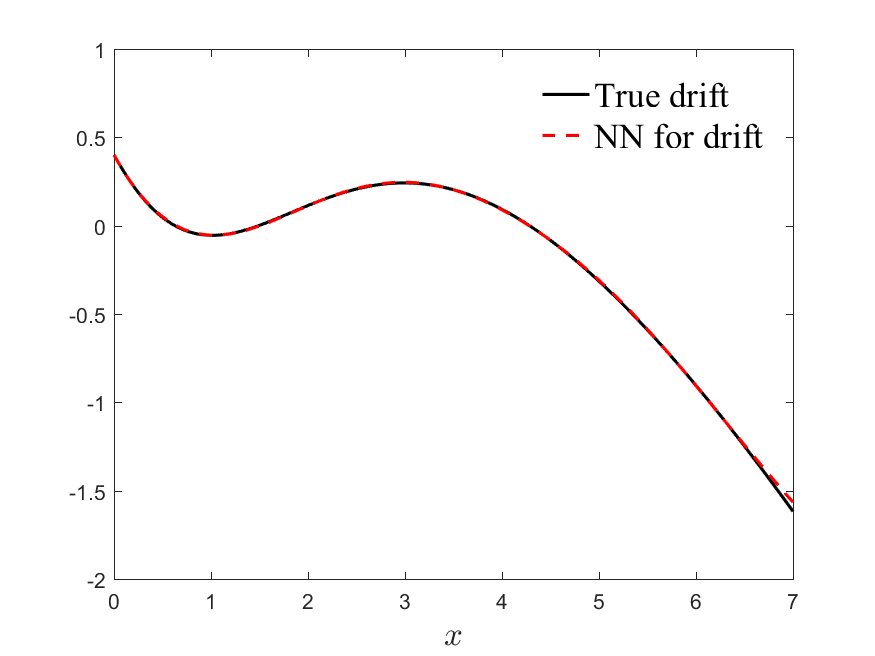}}
\end{minipage}
\hfill
\begin{minipage}[]{0.3 \textwidth}
 \leftline{\tiny\textbf{(b)}}
\centerline{\includegraphics[width=4.5cm]{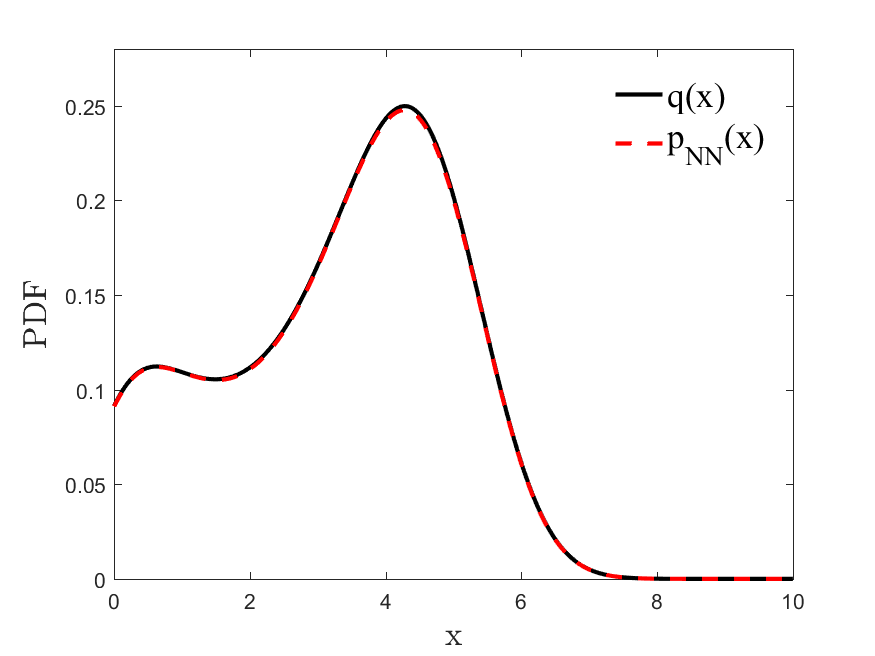}}
\end{minipage}
\hfill
\begin{minipage}[]{0.3 \textwidth}
 \leftline{\tiny\textbf{(c)}}
\centerline{\includegraphics[width=4.5cm]{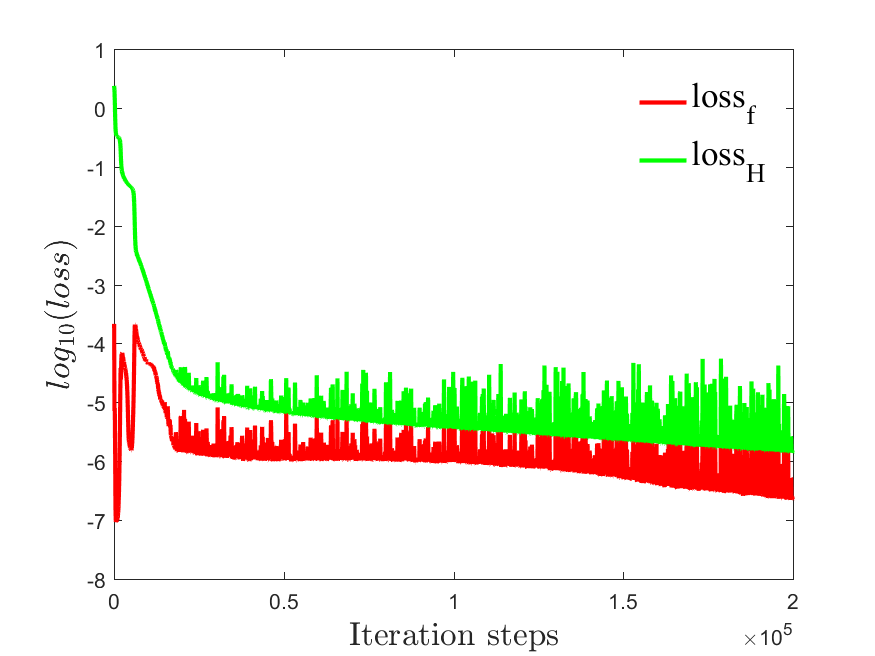}}
\end{minipage}
\caption{\textbf{ Unknown drift of Example 4:} (a) The learned drift term; (b) The learned probability; (c) The loss function.}
\label{case_bio}
\end{figure}

\begin{example}
This is a stochastic model for a transcription factor (i.e., a protein) concentration evolution in a certain gene regulation network
\[  dX=b(X)dt +  \sigma dB_t,  \]
with drift function $b(x)=\frac{k_f x^2}{x^2 + K_d} - k_d x + R_{bas}$.  Here the parameters are $K_d = 10$, $k_d = 1 min^{-1}$,  $k_f = 6 min^{-1}$, and $R_{bas} = 0.4 min^{-1}$.
We will find a drift function $b(x)$ and diffusion $\sigma$ so that the Hellinger distance $I(b(x),\sigma) = \frac12 \int_{\R} [\sqrt{p(x)} -\sqrt{q(x)}]^2 dx$ is minimized.
\end{example}
Given the noise intensity $\sigma=1$, two neural networks are used to approximate the drift term and stationary probability density respectively. Here we still choose $N_{H}=1001$, $N_f=10000$. The result of the drift is shown in Figure \ref{case_bio} (a), where the black line is the true drift term and the red line is the neural network result. We find that he learned drift term can fit the true result very well. And the learned probability density function is shown in Figure \ref{case_bio} (b), which fits very well with the observation probability $q(x)$. We also show the loss function of $loss_H$ and $loss_f$ in Figure \ref{case_bio} (c). We see that the loss function decreases fast with the iteration steps increasing.

\begin{figure}[h]
\begin{minipage}[]{0.3 \textwidth}
 \leftline{\tiny\textbf{(a)}}
\centerline{\includegraphics[width=4.5cm]{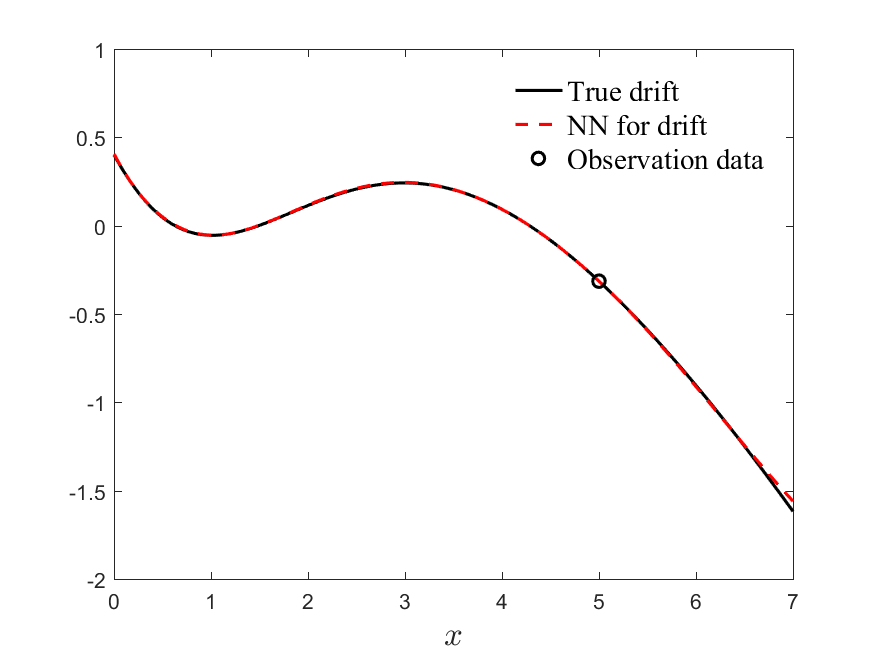}}
\end{minipage}
\hfill
\begin{minipage}[]{0.3 \textwidth}
 \leftline{\tiny\textbf{(b)}}
\centerline{\includegraphics[width=4.5cm]{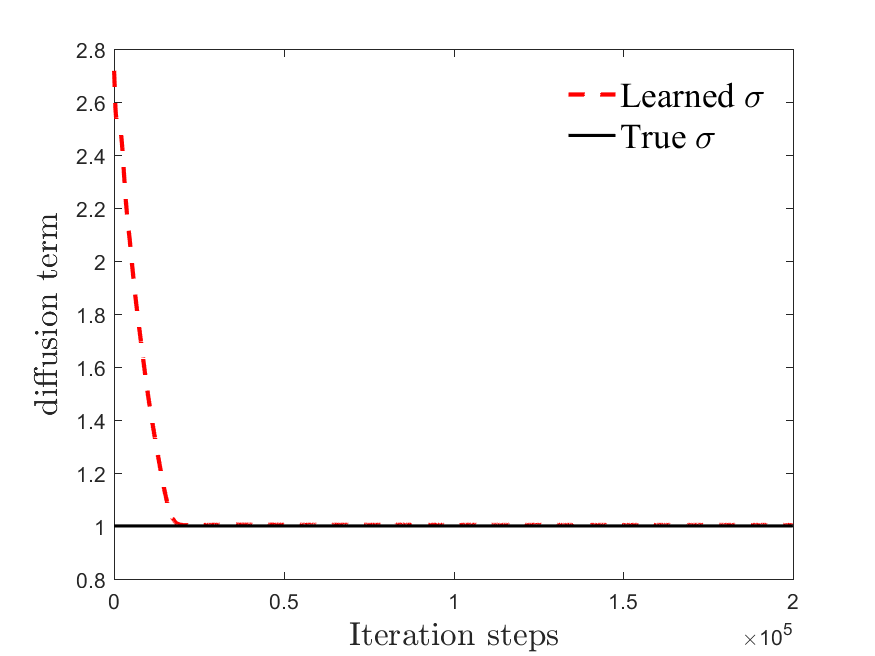}}
\end{minipage}
\hfill
\begin{minipage}[]{0.3 \textwidth}
 \leftline{\tiny\textbf{(c)}}
\centerline{\includegraphics[width=4.5cm]{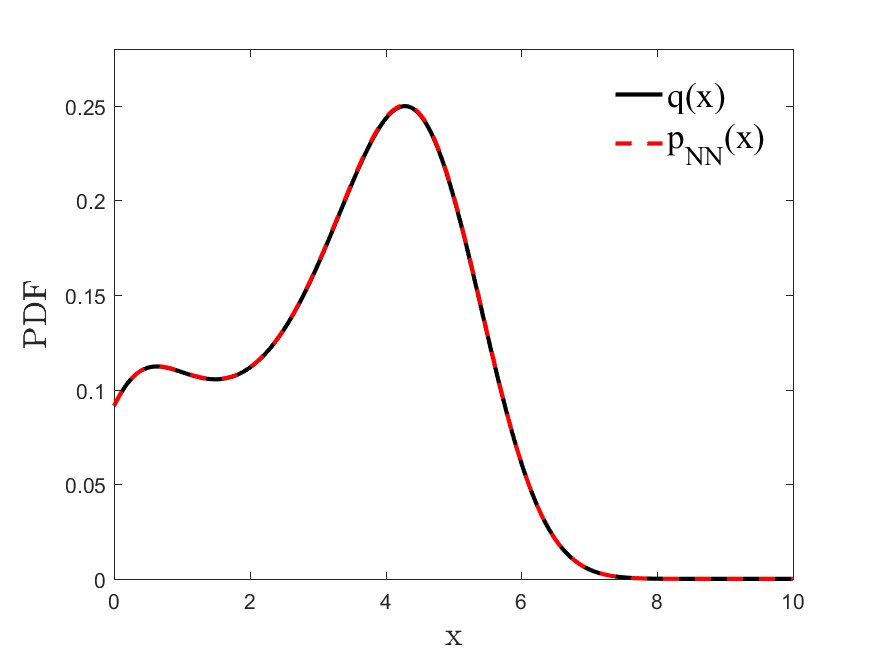}}
\end{minipage}
\caption{\textbf{ Unknown drift of Example 4:} (a) learned drift term with 1 observation data at $x=5$; (b) learned diffusion term; (c) learned probability density function.}
\label{case_bio_diffusion}
\end{figure}

For the unknown drift term and diffusion term case, we train the loss function \eqref{loss_case2} to get the optimal result. Only one observation data of drift term at $x=5$ is given. The results are shown in Figure \ref{case_bio_diffusion}.
 We present learned drift result in Figure \ref{case_bio_diffusion} (a), and the diffusion term evolution predictions as the iteration of the optimiser progresses in Figure \ref{case_bio_diffusion} (b). The neural network result of probability density is shown in Figure \ref{case_bio_diffusion} (c). From the figures,
we see that the drift and diffusion terms can be learned well.

\begin{example}
We now consider the following three dimensional stochastic dynamical systems with non-polynomial drift   \cite{chen2021solving}:
\begin{equation}
d\left( \begin{array}{ccc}
X_t\\
Y_t\\
Z_t
\end{array}
\right )=
\left( \begin{array}{c}
-\partial_{X_t} \Phi(X_t, Y_t,Z_t)\\
-\partial_{Y_t} \Phi(X_t, Y_t,Z_t)\\
-\partial_{Z_t} \Phi(X_t, Y_t,Z_t)\\
\end{array}
\right )dt+\left[ \begin{array}{ccc}
\sigma_1 & 0&0 \\
0& \sigma_2&0\\
0&0& \sigma_3\\
\end{array}
\right ]  d \left( \begin{array}{c}
B_{1,t}\\
B_{2,t}\\
B_{3,t}
\end{array}
\right ),
\end{equation}
where the potential
$\Phi(x,y,z)=-\frac{1}{2} \log[ (2\exp( \lambda_{01}(x-\lambda_{11} )^2 +\lambda_{02} (y-\lambda_{12})^2+\lambda_{03}(z-\lambda_{13})^2 )+\exp( \lambda_{04}(x-\lambda_{14})^2 +\lambda_{05}(y-\lambda_{15})^2+\lambda_{06}(z-\lambda_{16})^2 )  ]$,
$\lambda_{0i}=-5,-2.5,-5,-1,-1,-1$,
$\lambda_{1i}=1,1,1,-2,-1,-1$
and $\sigma_j=1$, where $i=1,2,...,6$ and $j=1,2,3$.
The ``observation" of the stationary probability density is $q(x,y,z)=1/Z  \exp( -2\Phi(x,y,z) )$, where $Z$ is the normalization parameter such that the integral of $q(x,y,z)$ on domain $\mathrm{R}^3$ is equal to 1.
Find the parameters in drift term and diffusion term so that the Hellinger distance $I = \frac12 \int_{\R^3} [\sqrt{p(x,y,z)} -\sqrt{q(x,y,z)}]^2 dxdydz$ is minimized.
\end{example}

We use neural network to approximate the stationary probability density. And here we   choose $N_{H}=50000$, $N_f=5000$.

\begin{figure}[h]
\begin{minipage}[]{0.3 \textwidth}
 \leftline{\tiny\textbf{(a)}}
\centerline{\includegraphics[width=4.5cm]{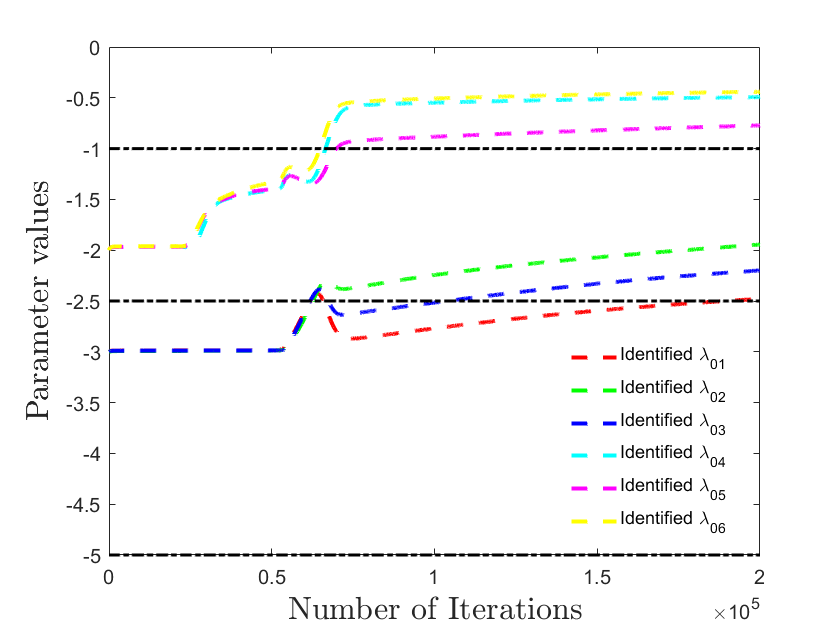}}
\end{minipage}
\hfill
\begin{minipage}[]{0.3 \textwidth}
 \leftline{\tiny\textbf{(b)}}
\centerline{\includegraphics[width=4.5cm]{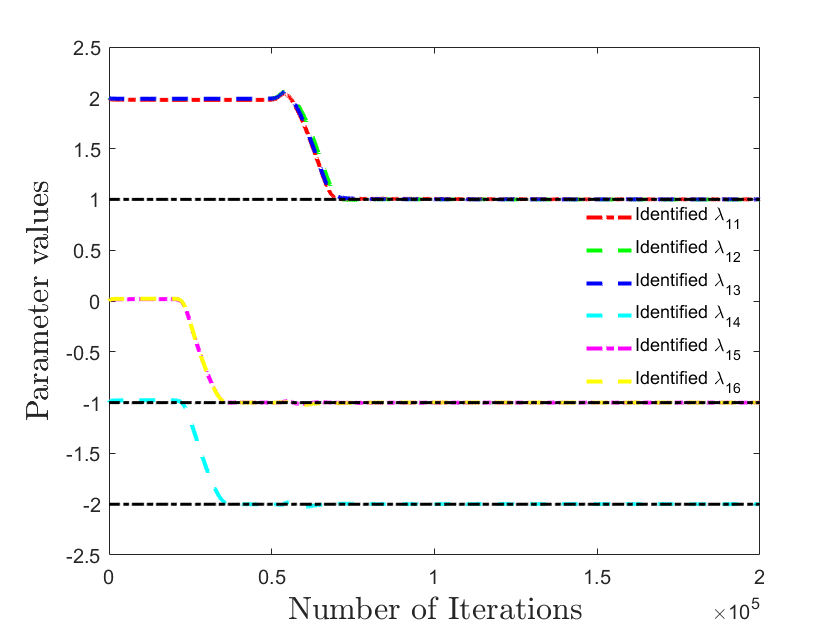}}
\end{minipage}
\hfill
\begin{minipage}[]{0.3 \textwidth}
 \leftline{\tiny\textbf{(c)}}
\centerline{\includegraphics[width=4.5cm]{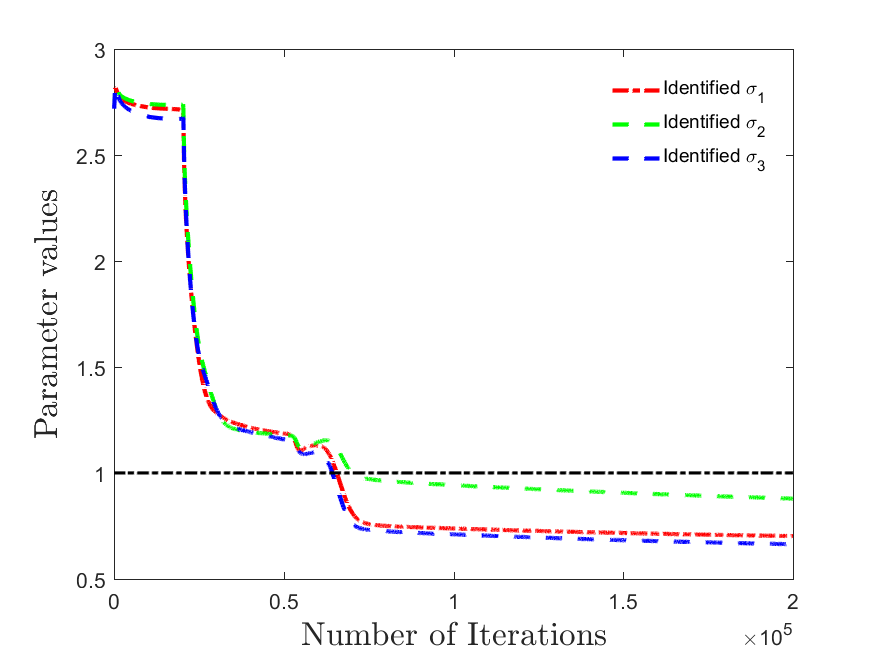}}
\end{minipage}
\caption{\textbf{Three dimensional results of Example 5:} (a) learned $\lambda_{0i}$; (b) learned $\lambda_{1i}$; (c) learned $\sigma_{j}$, where $i=1,2,...,6$ and $j=1,2,3$.}
\label{case_3d_all}
\end{figure}

First, we learn  all the parameters $\lambda_{0i}$, $\lambda_{1i}$ and $\sigma_j$, for $i=1,2,...6$ and $j=1,2,3$.
The results are shown in Figure \ref{case_3d_all}.
 In Figure \ref{case_3d_all}(a) and (c), the parameters of $\lambda_{0i}$ and drift term are learned not well. While the parameters $\lambda_{2i}$ can be learned well, see Figure \ref{case_3d_all} (b). So we will learn the parameter in the drift term and diffusion term respectively.

\begin{figure}[h]
\centerline{\includegraphics[width=8cm,height=6cm]{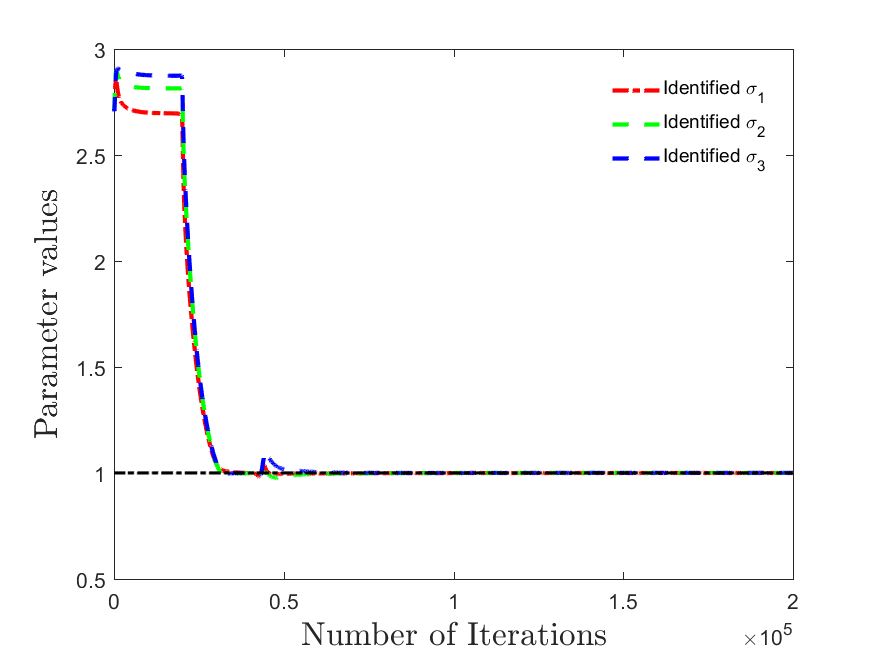}}
\caption{\textbf{ 3D result} the learned drift term of Example 5.}
\label{3d_diffion}
\end{figure}

On the one hand, we just learn the diffusion term given the drift term. The results are shown in Figure \ref{3d_diffion}. The parameters in the diffusion term approach to the true parameter as the number of iterations increases.

\begin{figure}[!htb]
\begin{minipage}[]{0.3 \textwidth}
 \leftline{\tiny\textbf{(Case I)}}
\centerline{\includegraphics[width=4.5cm]{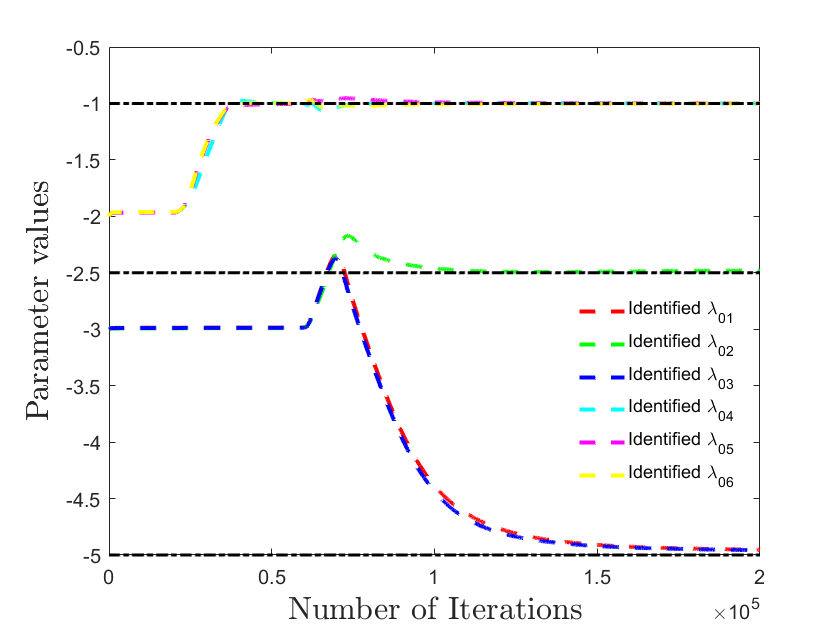}}
\end{minipage}
\hfill
\begin{minipage}[]{0.3 \textwidth}
 \leftline{\tiny\textbf{(Case II)}}
\centerline{\includegraphics[width=4.5cm]{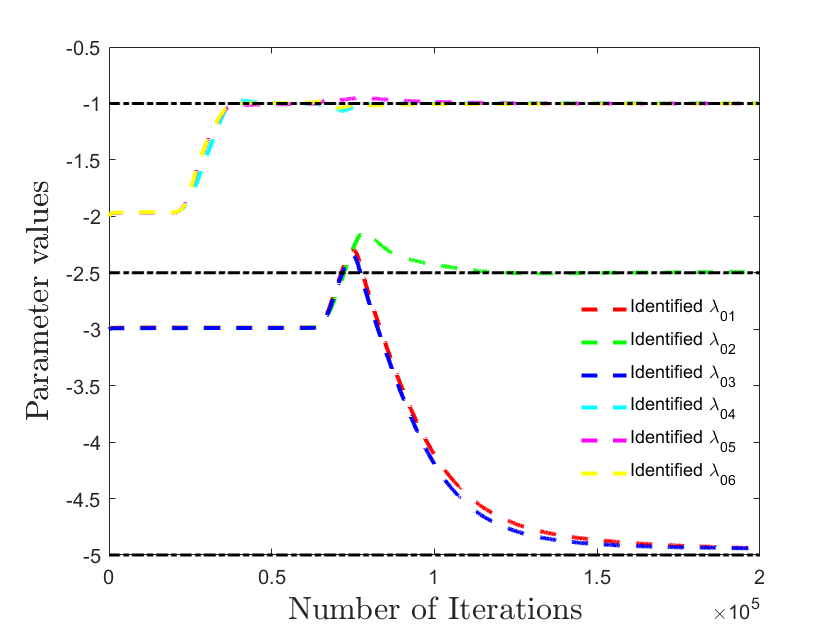}}
\end{minipage}
\hfill
\begin{minipage}[]{0.3 \textwidth}
 \leftline{\tiny\textbf{(Case III)}}
\centerline{\includegraphics[width=4.5cm]{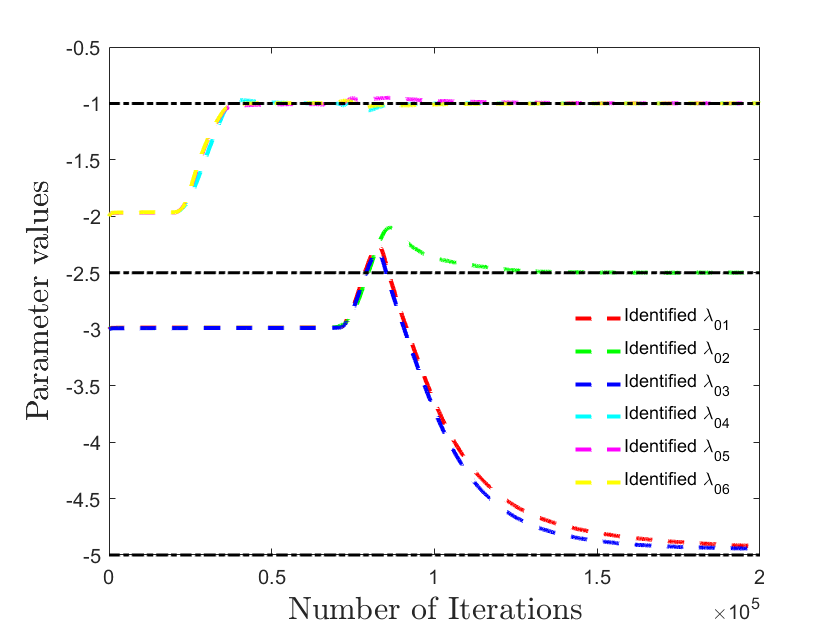}}
\end{minipage}

\begin{minipage}[]{0.3 \textwidth}
\centerline{\includegraphics[width=4.5cm]{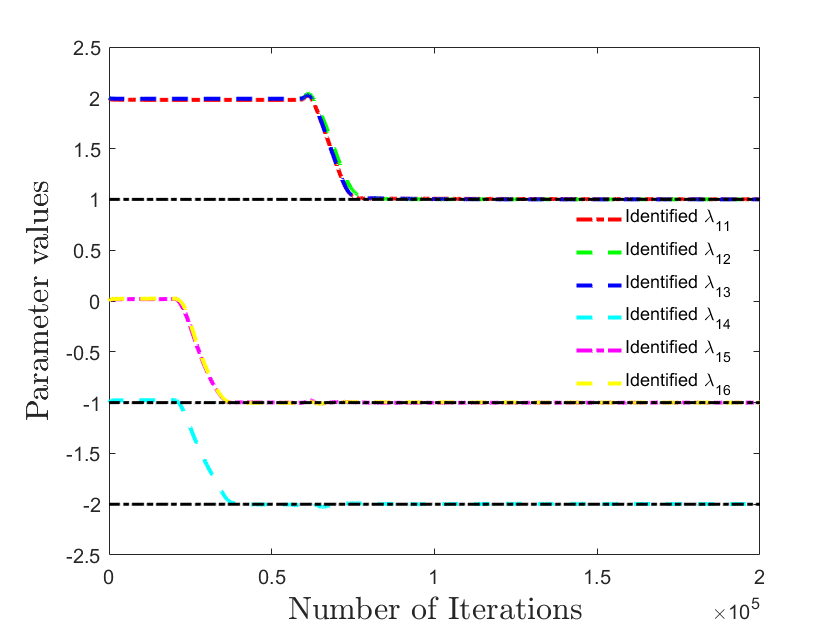}}
\end{minipage}
\hfill
\begin{minipage}[]{0.3 \textwidth}
\centerline{\includegraphics[width=4.5cm]{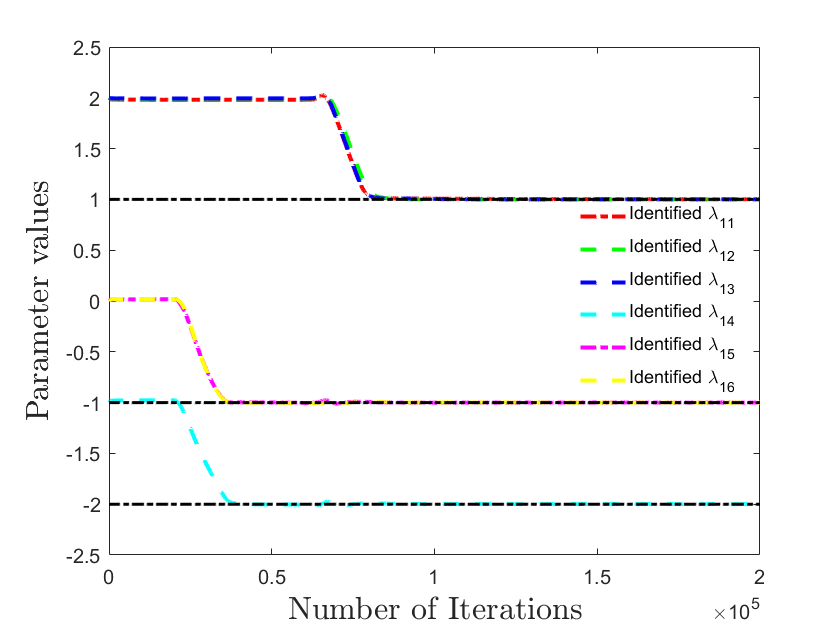}}
\end{minipage}
\hfill
\begin{minipage}[]{0.3 \textwidth}
\centerline{\includegraphics[width=4.5cm]{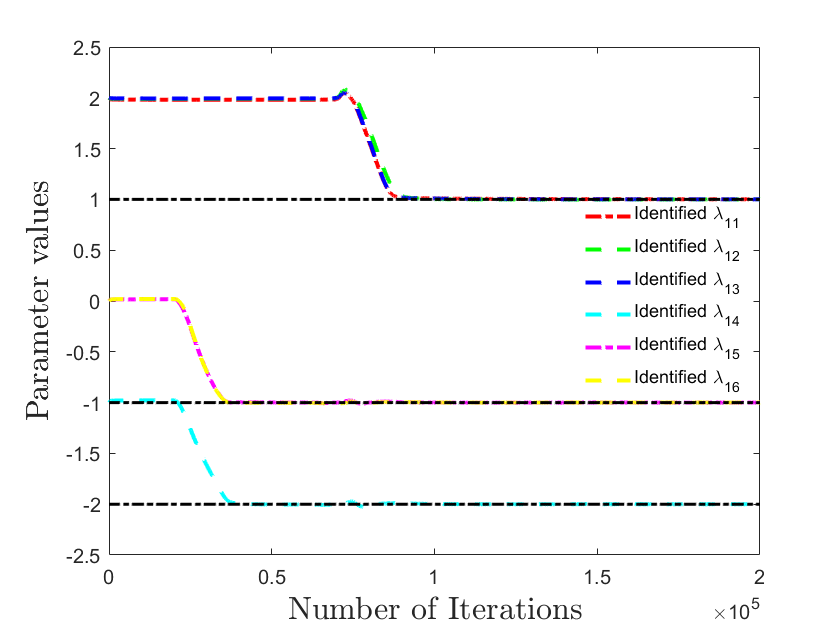}}
\end{minipage}
\caption{\textbf{Three dimensional results of Example 5.} Learn all parameters in the drift terms with perturbation. Case I: clean observation data of the PDF; Case II: $5\%$ noise observation data of the PDF; Case III: $10\%$ noise observation data of the PDF. Left: learned $\lambda_{0i}$; right: learned $\lambda_{1i}$, where $i=1,2,...,6$ and $j=1,2,3$.}
\label{3d_drift}
\end{figure}

On the other hand, we learn the parameters $\lambda_{0i}$ and $\lambda_{1i}$ in the drift term given the diffusion term. The results are shown in Figure \ref{3d_drift}. We learn the results for three cases. For case I: the observation data is clean, i.e. $q(x,y,z)$. For case II: the observation data is given with $5\%$ noise, i.e. $q(x,y,z)*(1+0.05\mathrm{N}(0,I))$, and for case III: the observation data is given with $10\%$ noise, i.e. $q(x,y,z)*(1+0.1\mathrm{N}(0,I))$. Here $\mathrm{N}(0,I)$ mean the standard normal distribution. The parameters we learned are well even the observation data has $10\%$ noise.

\begin{figure}[h]
\begin{minipage}[]{0.48 \textwidth}
\centerline{\includegraphics[width=6cm]{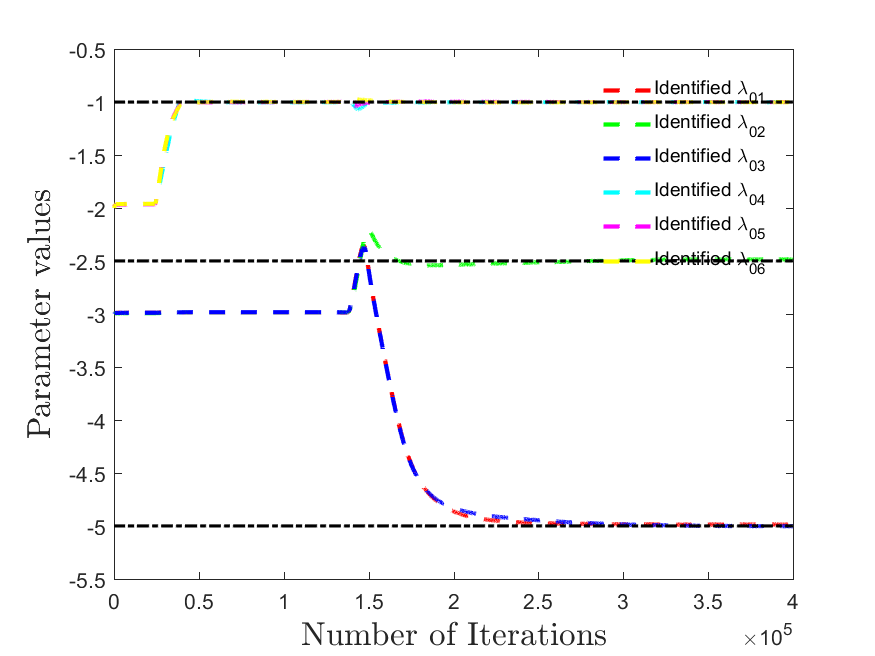}}
\end{minipage}
\hfill
\begin{minipage}[]{0.48 \textwidth}
\centerline{\includegraphics[width=6cm]{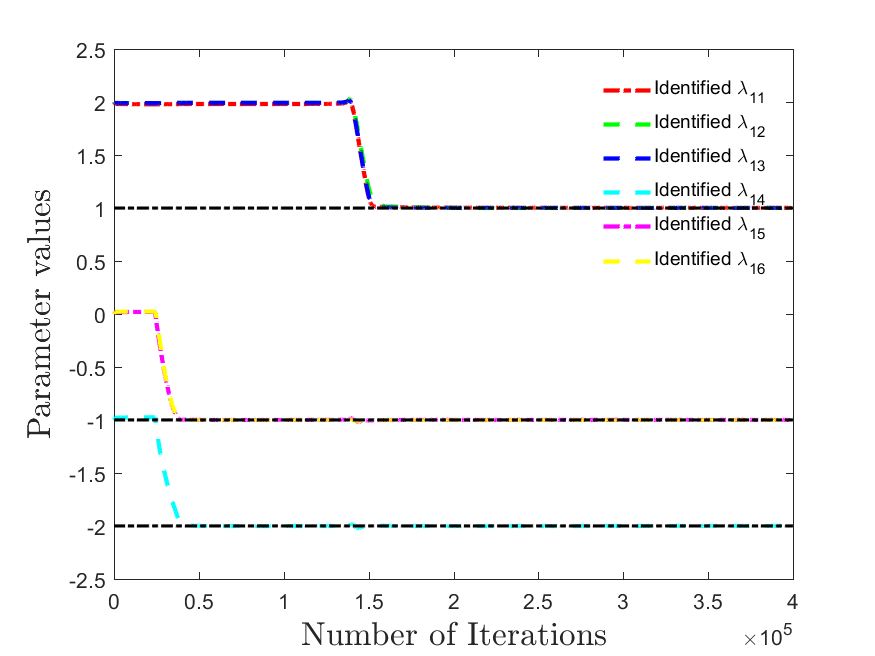}}
\end{minipage}
\caption{\textbf{Three dimensional results of Example 5.} Learn parameters in the drift terms using Jensen-Shannon distance. Left: learned $\lambda_{0i}$; right: learned $\lambda_{1i}$, where $i=1,2,...,6$ and $j=1,2,3$.}
\label{3d_drift_JS}
\end{figure}

In the following, We change the Hellinger distance to the Jensen-Shannon divergence in the loss function. The results of learned parameters in the drift term are shown in Figure \eqref{3d_drift_JS}. The unknown parameters can be learned well. While compared with the Hellinger distance, this method needs more iteration steps to train.

\begin{figure}[h]
\begin{minipage}[]{0.48 \textwidth}
\centerline{\includegraphics[width=6cm]{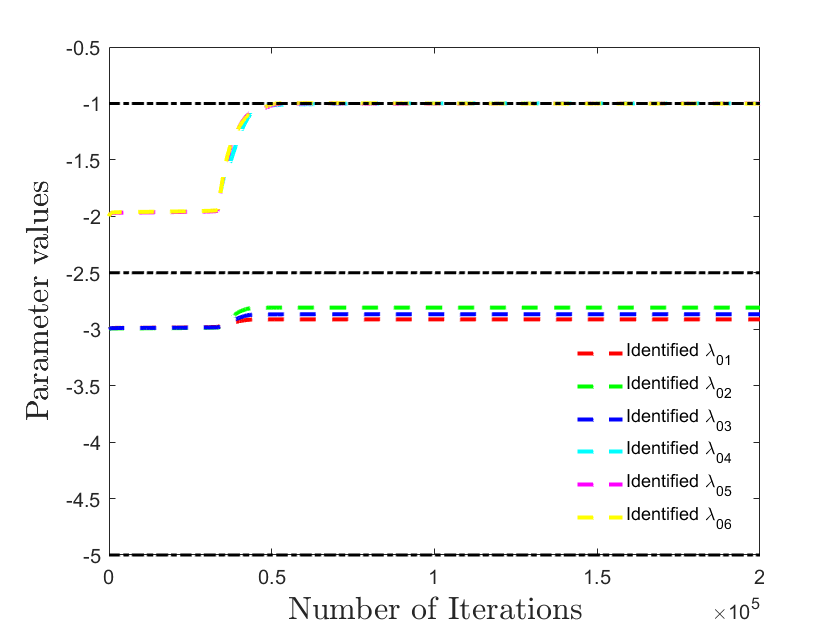}}
\end{minipage}
\hfill
\begin{minipage}[]{0.48 \textwidth}
\centerline{\includegraphics[width=6cm]{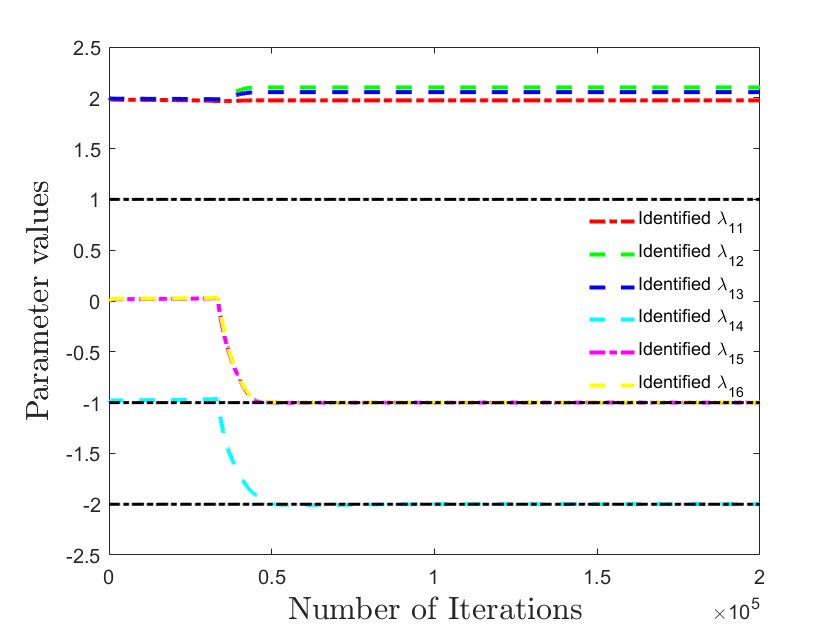}}
\end{minipage}
\caption{\textbf{Three dimensional results with PINN loss of Example 5.} Learn all the parameters in the drift term with clean observation data of the PDF. Left: learned $\lambda_{0i}$; right: learned $\lambda_{1i}$, where $i=1,2,...,6$ and $j=1,2,3$.}
\label{3d_drift_PINN}
\end{figure}

\begin{figure}[h]
\begin{minipage}[]{0.2 \textwidth}
\leftline{\tiny\textbf{(a1)}}
\centerline{\includegraphics[width=4.5cm]{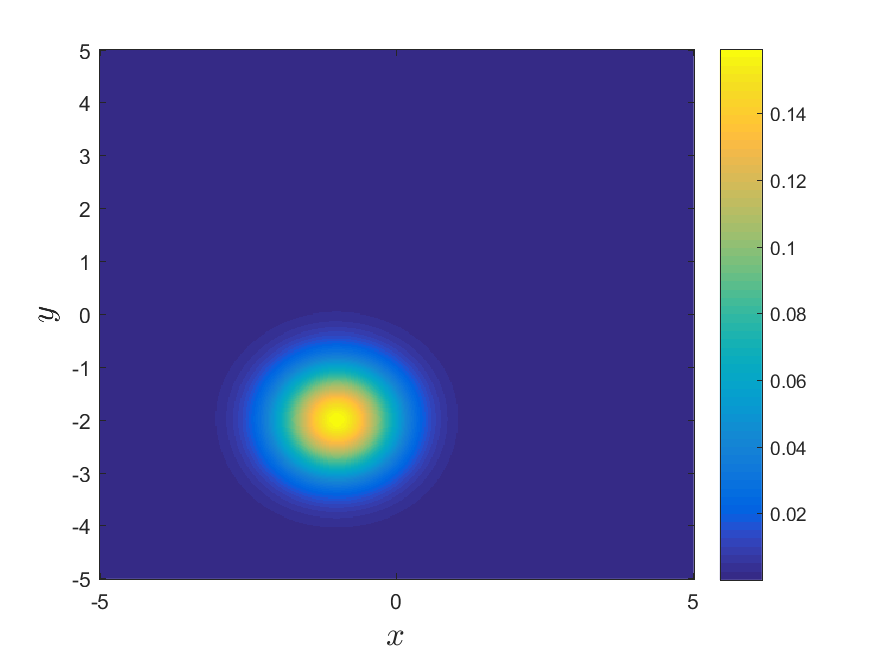}}
\end{minipage}
\hfill
\begin{minipage}[]{0.2 \textwidth}
\leftline{\tiny\textbf{(b2)}}
\centerline{\includegraphics[width=4.5cm]{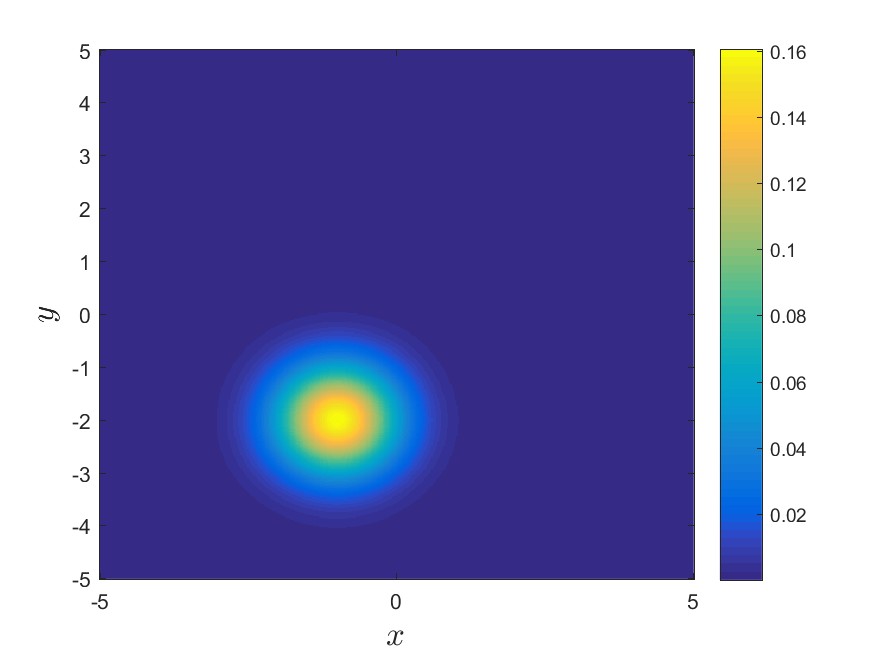}}
\end{minipage}
\hfill
\begin{minipage}[]{0.2 \textwidth}
\leftline{\tiny\textbf{(c1)}}
\centerline{\includegraphics[width=4.5cm]{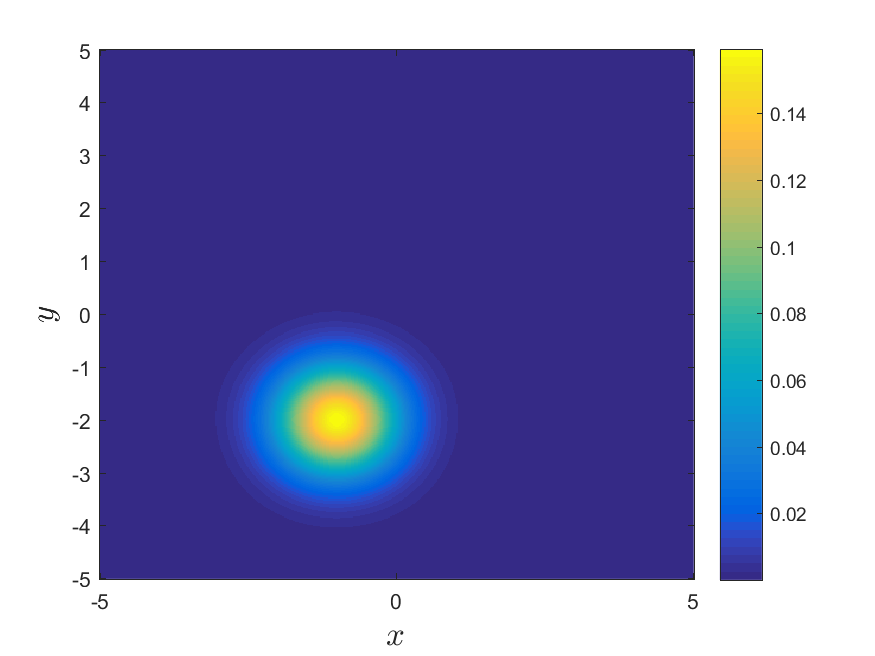}}
\end{minipage}
\hfill
\begin{minipage}[]{0.2\textwidth}
\leftline{\tiny\textbf{(d1)}}
\centerline{\includegraphics[width=4.5cm]{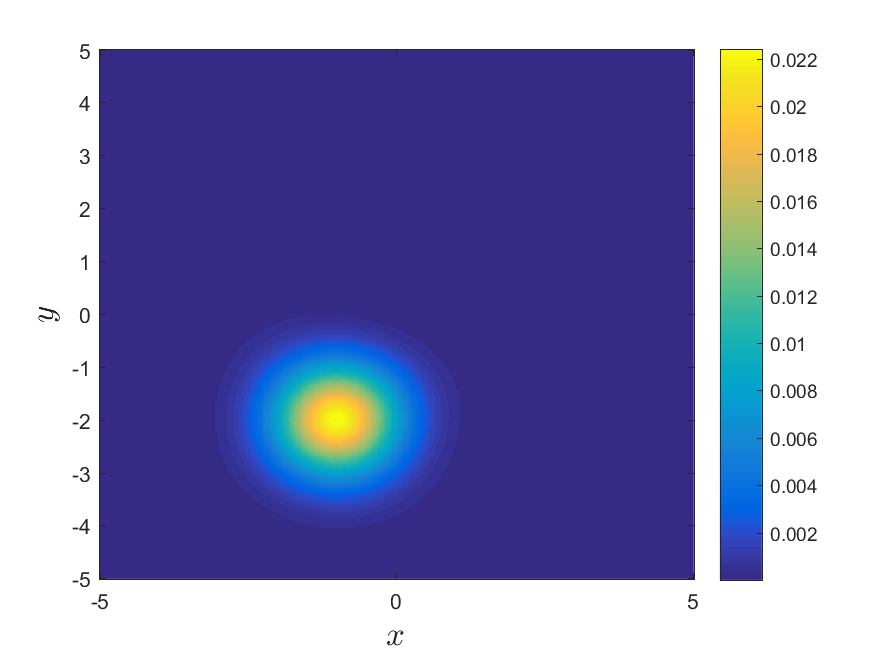}}
\end{minipage}
\begin{minipage}[]{0.2 \textwidth}
\leftline{\tiny\textbf{(a2)}}
\centerline{\includegraphics[width=4.5cm]{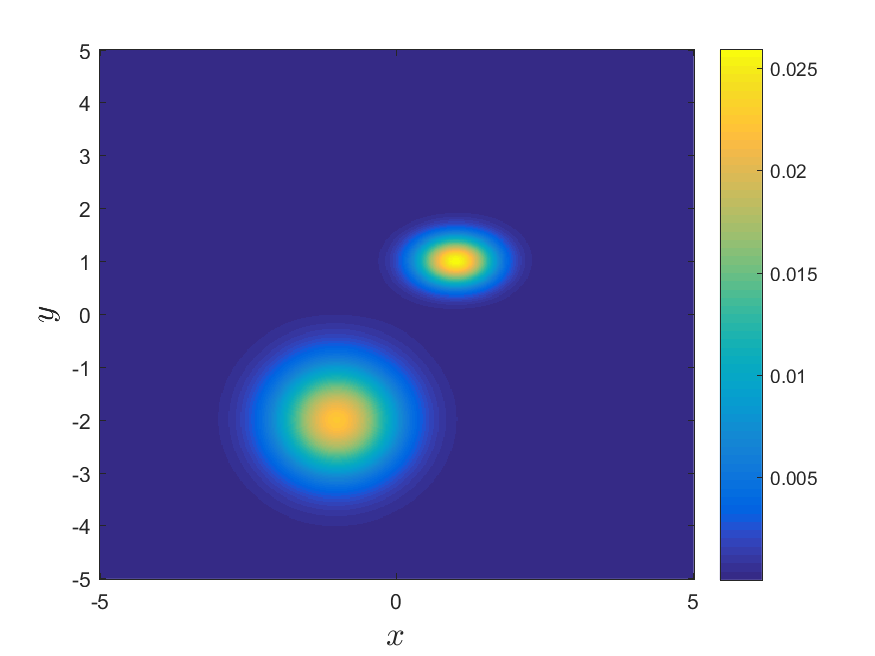}}
\end{minipage}
\hfill
\begin{minipage}[]{0.2 \textwidth}
\leftline{\tiny\textbf{(b2)}}
\centerline{\includegraphics[width=4.5cm]{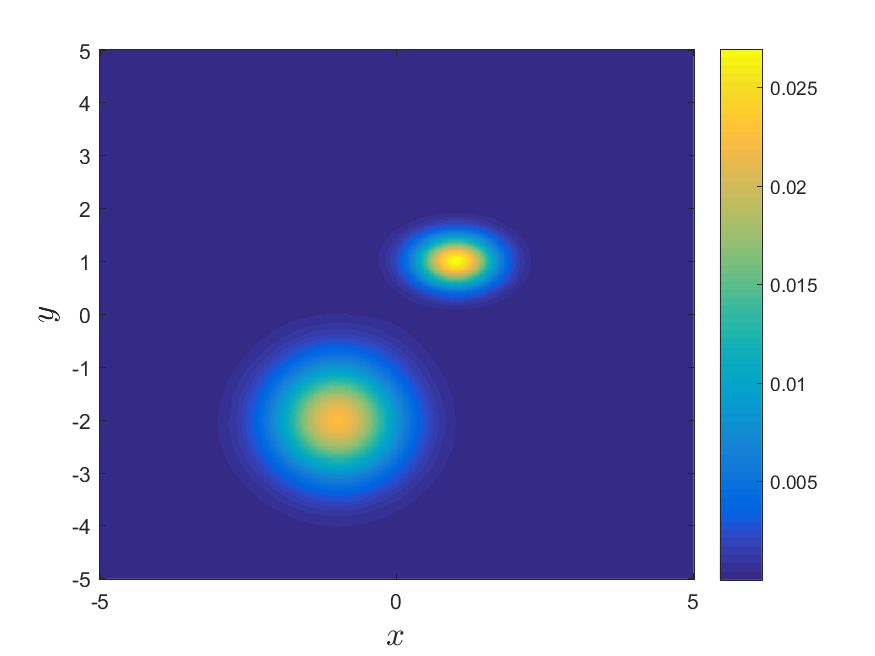}}
\end{minipage}
\hfill
\begin{minipage}[]{0.2 \textwidth}
\leftline{\tiny\textbf{(c2)}}
\centerline{\includegraphics[width=4.5cm]{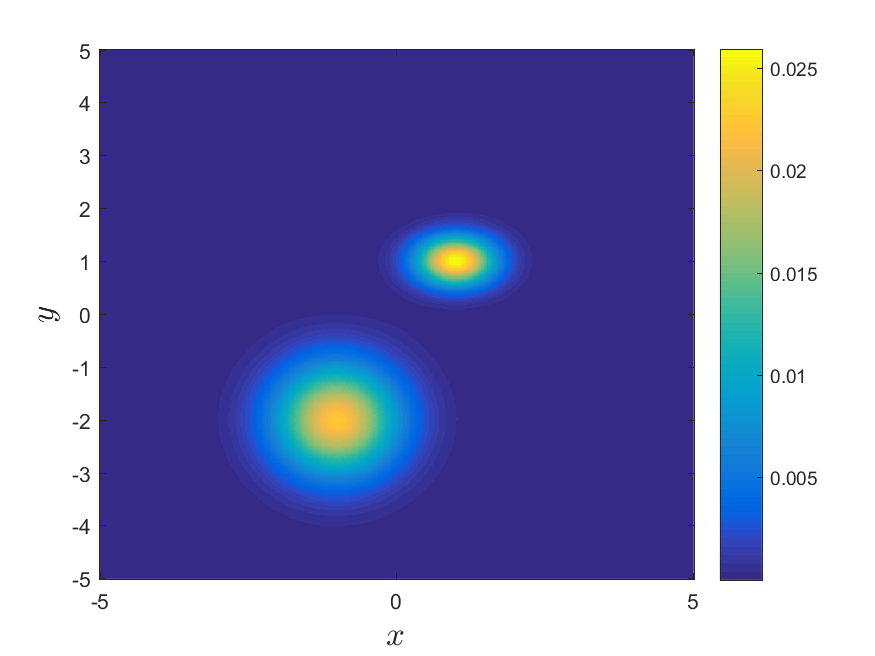}}
\end{minipage}
\hfill
\begin{minipage}[]{0.2\textwidth}
\leftline{\tiny\textbf{(d2)}}
\centerline{\includegraphics[width=4.5cm]{pdf_NN_PINN_015.png}}
\end{minipage}
\begin{minipage}[]{0.2 \textwidth}
\leftline{\tiny\textbf{(a3)}}
\centerline{\includegraphics[width=4.5cm]{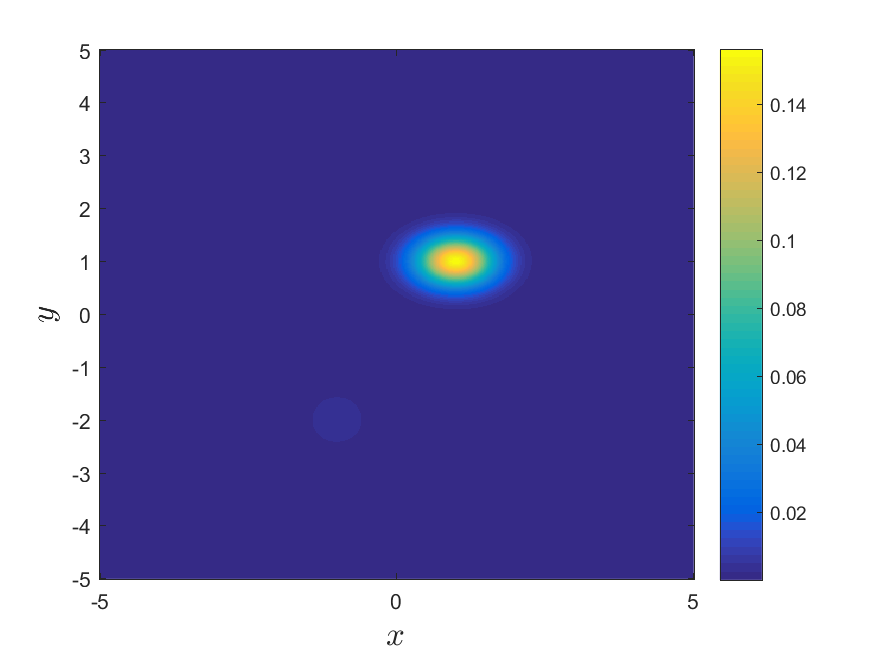}}
\end{minipage}
\hfill
\begin{minipage}[]{0.2 \textwidth}
 \leftline{\tiny\textbf{(b3)}}
\centerline{\includegraphics[width=4.5cm]{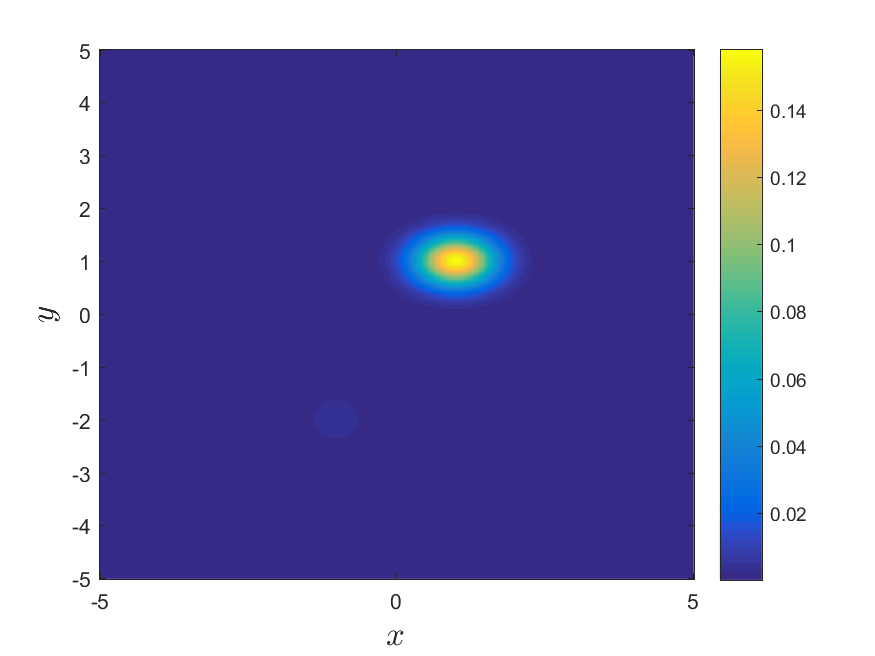}}
\end{minipage}
\hfill
\begin{minipage}[]{0.2 \textwidth}
\leftline{\tiny\textbf{(c3)}}
\centerline{\includegraphics[width=4.5cm]{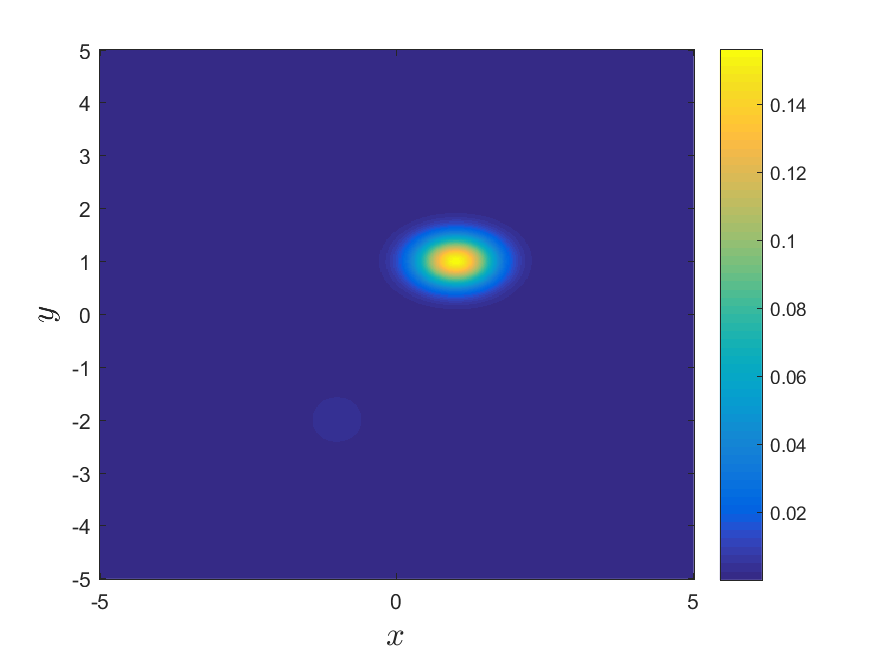}}
\end{minipage}
\hfill
\begin{minipage}[]{0.2\textwidth}
\leftline{\tiny\textbf{(d3)}}
\centerline{\includegraphics[width=4.5cm]{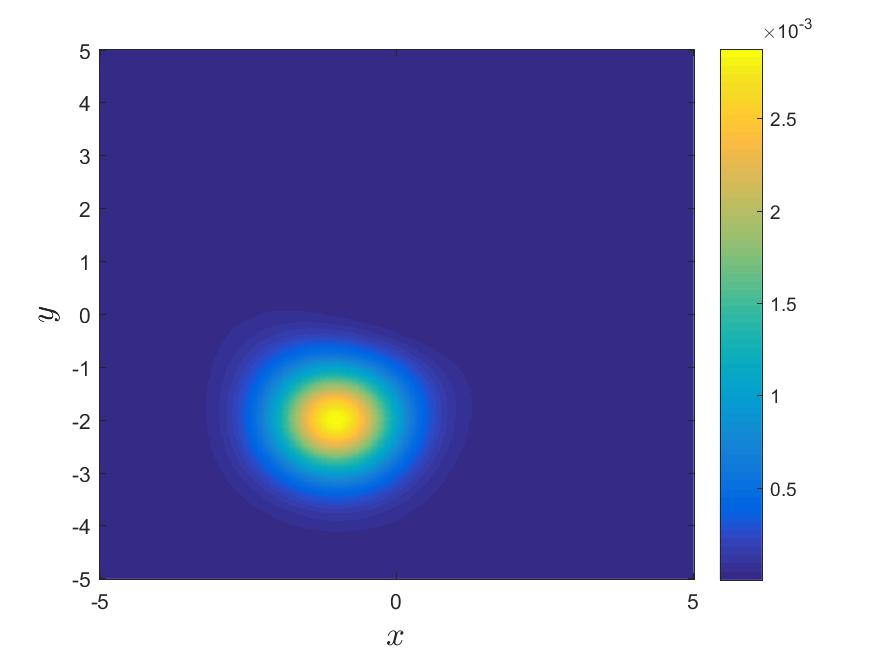}}
\end{minipage}
\caption{\textbf{The probability density function of Example 5.} (a1)-(a3): the true PDF with $z=-1,0.5,1$; (b1)-(b3): the learned PDF using Hellinger distance (\ref{loss}); (c1)-(c3): the learned PDF using Jensen-Shannon divergence loss (\ref{loss_JS} ); (d1)-(d3): the learned PDF using PINN loss (\ref{loss_PINN}).  }
\label{3d_PDF}
\end{figure}

Here we also compare our results with the traditional physics informed neural network (PINN) with the case of learning the parameter in the drift term.
The results are shown in Figure \ref{3d_drift_PINN}. Only several parameters in the drift term can be learned well using PINN method. Compared with PINN method using mean square error \eqref{loss_PINN}, our loss with Hellinger distance \eqref{loss} would get better results. We use the neural network to approximate the probability density function and plot the learned probability density function when $z=-1,0.5,1$ using different methods. The results are shown in Figure \ref{3d_PDF}. Compared with the true probability density function, the proposed method with Hillinger distance and Jensen-Shannon divergence works better than the mean square distance.

\begin{figure}[!h]
\centerline{\includegraphics[width=16cm]{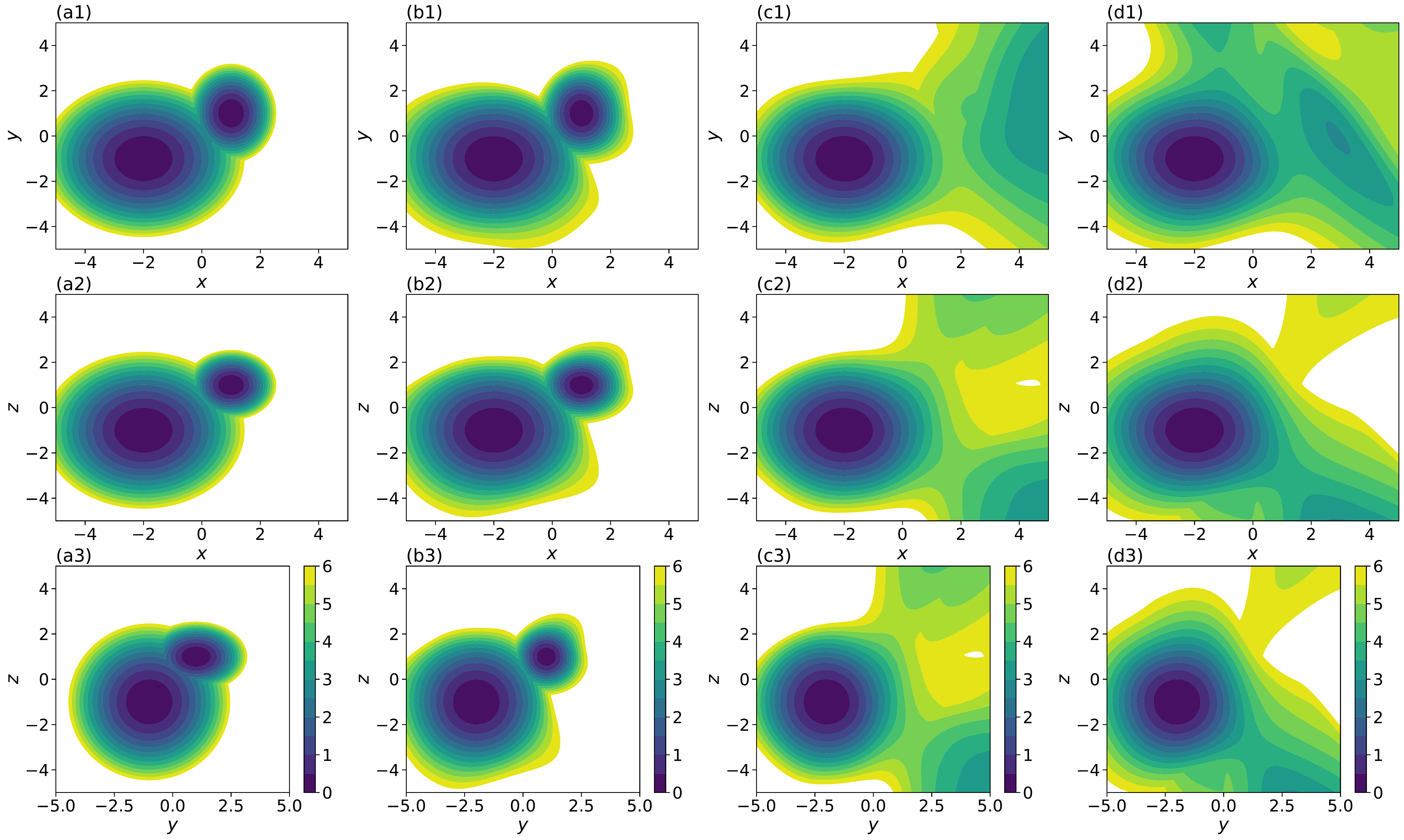}}
\caption{\textbf{The learned potential of Example 5.} (a1)-(a3): the true potential; (b1)-(b3): the learned potential using Hellinger distance (\ref{loss}); (c1)-(c3): the learned potential using Jensen-Shannon divergence loss (\ref{loss_JS}); (d1)-(d3): the learned potential using PINN loss (\ref{loss_PINN}).  }
\label{3d_potential}
\end{figure}

For the unknown drift term, we   use our proposed method to recover the drift term. The results are shown in Figure \ref{3d_potential}. The first row of the Figure  \ref{3d_potential} (a1,b1,c1,d1) are the the projection of $\Phi(x,y,z)$ var minimization of $z$, i.e. $\min_z \Phi(x,y,z)$. The second and third rows are projected on $(x,z)$ domain and $(y,z)$ domain separately. The true projections are shown in Figure \ref{3d_potential}(a). The results with Hellinger distance are shown in Figure \ref{3d_potential}(b). Comparing with the true potential, this example illustrates our method with Hellinger distance works well. However, it is difficult to recover the SDE with the Jensen-Shannon divergence (see Figure \ref{3d_potential}(c)) and PINN losses (see Figure \ref{3d_potential}(d)). This indicates the efficiency of our proposed method using the Hellinger distance in this example.

\begin{figure}[h]
\begin{minipage}[]{0.48 \textwidth}
\centerline{\includegraphics[width=6cm]{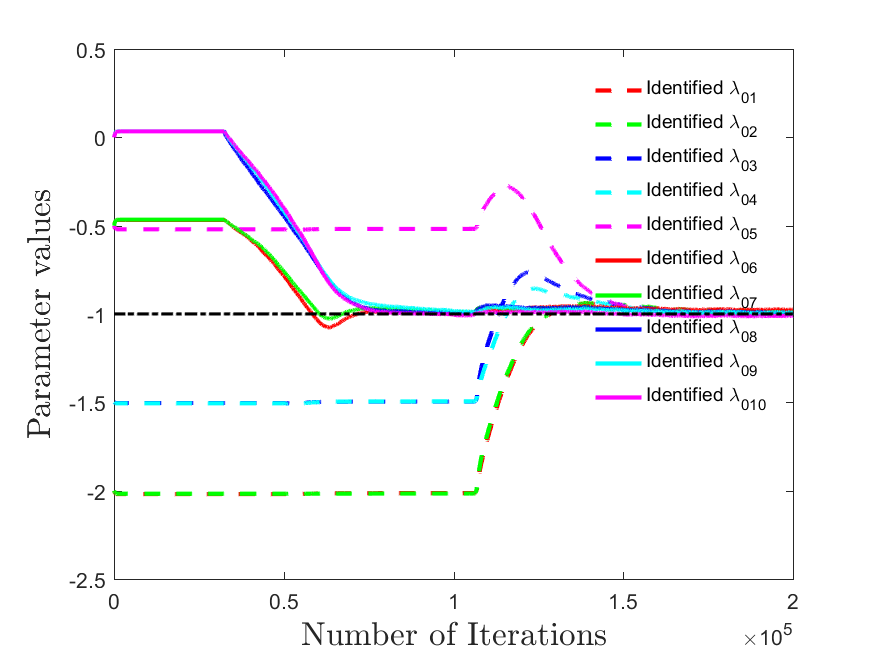}}
\end{minipage}
\hfill
\begin{minipage}[]{0.48 \textwidth}
\centerline{\includegraphics[width=6cm]{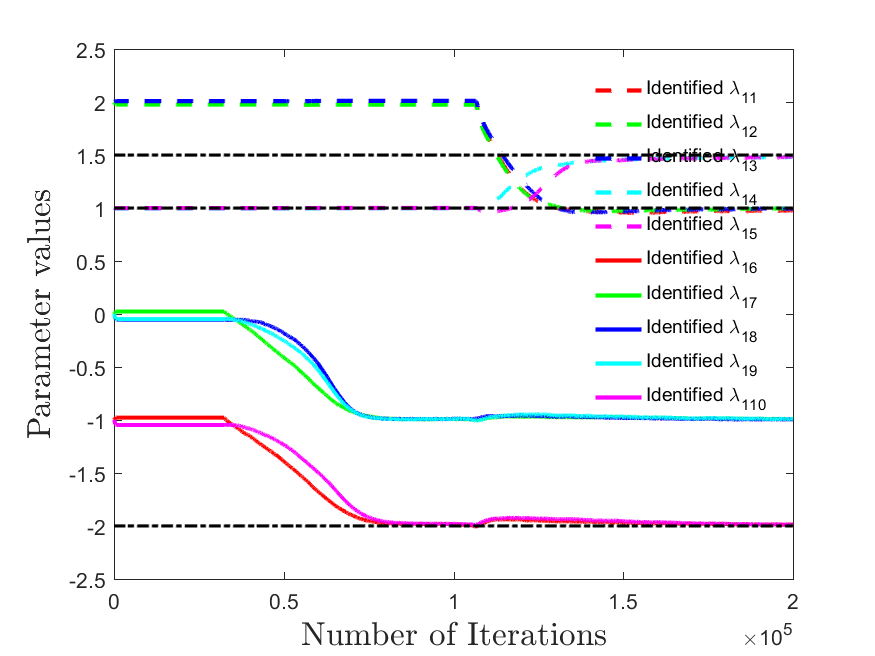}}
\end{minipage}
\caption{\textbf{Five dimensional results of Example 6.} Learn all the parameters in the drift term with clean observation data of the PDF. Left: learned $\lambda_{0i}$; right: learned $\lambda_{1i}$, where $i=1,2,...,10$.}
\label{5d_drift}
\end{figure}

\begin{example}
Finally, We now consider the following five dimensional stochastic dynamical systems with non-polynomial drift:
\begin{equation} \nonumber
d\left( \begin{array}{ccccc}
X_t\\
Y_t\\
Z_t\\
V_t\\
W_t
\end{array}
\right )=
\left( \begin{array}{c}
-\partial_{X_t} \Phi(X_t, Y_t, Z_t, V_t, W_t)\\
-\partial_{Y_t} \Phi(X_t, Y_t, Z_t, V_t, W_t)\\
-\partial_{Z_t} \Phi(X_t, Y_t, Z_t, V_t, W_t)\\
-\partial_{V_t} \Phi(X_t, Y_t, Z_t, V_t, W_t)\\
-\partial_{W_t} \Phi(X_t, Y_t, Z_t, V_t, W_t)\\
\end{array}
\right )dt+\left[ \begin{array}{ccccc}
\sigma_1 & 0&0 &0&0 \\
0& \sigma_2&0&0&0\\
0&0& \sigma_3&0&0\\
0&0&0 &\sigma_4&0\\
0&0&0 &0&\sigma_5\\
\end{array}
\right ]  d \left( \begin{array}{c}
B_{1,t}\\
B_{2,t}\\
B_{3,t}\\
B_{4,t}\\
B_{5,t}
\end{array}
\right ),
\end{equation}

where the potential
$\Phi(x,y,z, v, w)=-\frac{1}{2} \log[ exp( (\lambda_{01}(x-\lambda_{11}) +\lambda_{02}(y-\lambda_{12})+\lambda_{03}(z-\lambda_{13}) +\lambda_{04}(v-\lambda_{14})+\lambda_{05}(w-\lambda_{15})))+exp( (\lambda_{06}(x-\lambda_{16})+\lambda_{07}(y-\lambda_{17}))+\lambda_{08}(z-\lambda_{18}) +\lambda_{09}(v-\lambda_{19})+\lambda_{10}(w-\lambda_{110}))   ]$,
$\lambda_{0i}=-1$,
$\lambda_{1i}=(1,1,1,1.5,1.5,-2,-1,-1,-1,-2) $
 where $i=1,2,\cdots,10$,  and $\sigma_j=1$, $j=1,2,3,4,5$.
The ``observation" of the stationary probability density is $q(x,y,z,v,w)=1/Z  \exp( -2\Phi(x,y,z,v,w) )$, where $Z$ is the normalization parameter such that the integral of $q(x,y,z,v,w)$ on domain $\R^5$ is equal to 1.
Find the parameters in drift term so that the Hellinger distance $I = \frac12 \int_{\R^5} [\sqrt{p(x,y,z,v,w)} -\sqrt{q(x,y,z,v,w)}]^2 dxdydz$ is minimized.
\end{example}
We use a neural network to approximate the stationary probability density. And here we choose $N_{H}=50000$, $N_f=5000$.
We learn the parameters $\lambda_{0i}$ and $\lambda_{1i}$ with $i=1,2,\cdots,10$ in the drift term given the diffusion term. The results are shown in Figure \ref{5d_drift},  indicating that our method also works for five dimensional case.

Remark: With only one  observation  trajectory   data, we first use kernel density estimation to approximate the stationary probability density function,
and then learn the SDE model. This works well for one dimension,  while for  high dimensional cases, we need   data on multiple trajectories.  


\section{Discussion}\label{sec:4}
Based on minimizing Hellinger distance between two probability distributions, we have devised a data-driven method to extract stochastic dynamical systems models from observation data of either
long time trajectories or   stationary probability distributions.
Our numerical results in one, three dimensionaland five dimensional  examples have verified that this method is feasible. We may also take other distances, in the space of probability distributions, in our method.  Indeed, we have also tried  our method using  Jensen-Shannon divergence and mean square distance.


 In principle, we may extend our method to learn high dimensional stochastic dynamical systems.
But when dealing with a high dimensional case,  the larger search space for the Fokker-Planck equation makes it difficult to train the neural network. We will try to use parallel computing \cite{parallel-PINN} or active sampling \cite{Active_learning} to train the neural network.
Moreover, our method is more stringent on data requirements for higher dimensions cases.
Therefore, in the future, we are going to explore other method to recover stochastic differential equation models.
We will also try to use our method to learn the stochastic differential equations  with non-Gaussian L\'{e}vy noise.

This approach leads to a data-driven stochastic dynamical systems study of random phenomena, as we can further examine dynamical behaviors of the leaned stochastic governing models \cite{Duan2015}.

\section*{Acknowledgements}
We would like to thank Xi Chen for helpful discussions. This work is supported by the National Natural Science Foundation of China (NSFC) (Grant No.11901536, 12141107). Xiaoli Chen is supported by the Ministry of Education, Singapore, under its Research Centre of Excellence award to the Institute for Functional Intelligent Materials (I-FIM, project No. EDUNC-33-18-279-V12)

\section*{CRediT authorship contribution statement}
X. Chen: Conceptualization, Software, Formal analysis, Writing-original draft. H. Wang: Conceptualization, Formal analysis,
Writing-original draft. J. Duan:  Discussion and Suggestion.

\section*{Declaration of competing interest}
 The authors declare that they have no known competing financial interests or personal relationships that could have appeared to influence the work reported in this paper.


\bibliography{references}

\begin{thebibliography}{100}
\bibitem{Arnold2003}
 L. Arnold,  {\em Random Dynamical Systems}. New York,  Springer,(2003) Corrected 2nd printing.

\bibitem{Pas}C. Pasquero, E.  Tziperman,   Statistical parameterization of heterogeneous oceanic convection, \emph{Journal of Physical Oceanography}. (2007) 37: 214-229.


\bibitem{Sura}C. Penland,  P.  Sura,  Sensitivity of an ocean model to``details" of stochastic forcing. In \emph{Proc. ECMWF Workshop on Represenation of Subscale Processes using Stochastic-Dynamic
Models}.  Reading, England, (2005)6-8 June.


\bibitem{GaoDuan}T. Gao,  J. Duan,  Quantifying model uncertainty in dynamical systems driven by non-gaussian L\'{e}vy stable noise with observations on mean exit time or escape probability. Communications in Nonlinear Science and Numerical
    Simulation. 39(2016) 1-6.


\bibitem{Wu}D. Wu, M. Fu,  J. Duan,  Discovering mean residence time and escape probability from data of stochastic dynamical systems. Chaos.  29(9)(2019)093122.
\bibitem{Hung}
C. L. Hung,   X. Zhang,  N. Gemelke,  C.  Chin,  Observation of scale invariance and universality in two-dimensional Bose gases. Nature.  470(7333)(2011) 236-239.

\bibitem{Hairapetian}
G. Hairapetian,   R. Stenzel,   Observation of a stationary, current-free double layer in a plasma. Physical Review Letters. 65(2)(1990)175.


\bibitem{Yarmchuk}
  E. J. Yarmchuk,  M. J. V. Gordon,    R. E. Packard, Observation of Stationary Vortex Arrays in Rotating Superfluid Helium. Physical Review Letters.  43(3)(1979)214-217.


\bibitem{Gefen}
O. Gefen,  O. Fridman, I. Ronin,    N. Q. Balaban, Direct observation of single stationary-phase bacteria reveals a surprisingly long period of constant protein production activity. Proceedings of the National Academy of Sciences.  111(1)(2014)556-561.

\bibitem{ArnoldWishtutz}
L. Arnold,  V.  Wishtutz,  Stationary solutions of linear systems with additive and multiplicative noise. Stochastics-an International Journal of Probability \& Stochastic Processes. 7(1-2)(1982)133-155.


\bibitem{Liberzon}D. Liberzon,  R. W.  Brockett,  Nonlinear feedback systems perturbed by noise: steady-state probability distributions and optimal control. Automatic Control IEEE Transactions on. 45(6)(2000)1116-1130.


\bibitem{Gray}A. H. Gray,  Uniqueness of Steady-State Solutions to the Fokker-Planck Equation. Journal of Mathematical Physics.  6(4)(1965)644-647.

\bibitem{Khasminskii}
 R.  Khasminskii,  Stochastic stability of differential equations,(2012)  Springer.

\bibitem{Schmalfuss}
B. Schmalfuss,   Lyapunov functions and non-trivial stationary solutions of stochastic differential equations. Dynamical Systems: An International Journal.  16(4)(2001) 303-317.

\bibitem{Gerber} F. Gerber,   D. W. Nychka, Fast covariance parameter estimation of spatial Gaussian process models using neural networks. \emph{Statistics and Probability}.(2021) 10(1).

\bibitem{Batz} P. Batz,  A. Ruttor,   M.  Opper,  Variational estimation of the drift for stochastic differential equations from the empirical density.  	Journal of Statistical Mechanics- Theory and Experiment.  (8)(2016) 083404.

\bibitem{batz2018approximate}
 P. Batz,  A. Ruttor,  M.  Opper,    Approximate Bayes learning of stochastic differential equations. Physical Review  E.  98(2)(2018) 022109.



\bibitem{opper2017estimator}
M. Opper,     An estimator for the relative entropy rate of path measures for stochastic differential equations. Journal of Computational Physics.  330(2017) 127-133.

\bibitem{opper2019variational}
M. Opper,   Variational inference for stochastic differential equations.  Annals of Physics.  531(3)(2019) 1800233.

\bibitem{ryder2018black}
T. Ryder,  A. Golightly,  A. S. McGough,   D. Prangle,   Black-box variational inference for stochastic differential equations. ICML.(2018)  4423-4432.

\bibitem{boninsegna2018sparse}
L. Boninsegna, F.  N\"{u}ske,  C.  Clementi,   Sparse learning of stochastic dynamical equations.  Journal of chemical physics.  148(24)(2018) 241723.

\bibitem{tabar2019analysis}
  R. Tabar,   Analysis and data-based reconstruction of complex nonlinear dynamical systems, Berlin/Heidelberg, (2019)Germany: Springer.



\bibitem{li2020scalable}
 X. Li,  T. K. L. Wong,  R. T. Chen,   D. Duvenaud,   Scalable gradients for stochastic differential equations. AISTATS. (2020) 3870-3882.



\bibitem{jia2019neural}
 J. Jia,    A. R. Benson,   Neural jump stochastic differential equations, NIPS.  (2019)32.



\bibitem{YangLiu-GAN-SODE}
L. Yang,  C. Daskalakis,  G. E. Karniadakis,   Generative Ensemble Regression: Learning Particle Dynamics from Observations of Ensembles with Physics-Informed Deep Generative Models. SIAM Journal on Scientific Computing.  44(1)(2022) B80-B99.

\bibitem{Zhang2020}
  H. Zhang,  Y. Xu,  Y. Li,   J. Kurths,   Statistical solution to SDEs with $\alpha$-stable L\'{e}vy noise via deep neural network. International Journal of Dynamics and Control. 8(4)(2020) 1129-1140.


\bibitem{Xu2020}
 Y. Xu,    H. Zhang,  Y. Li,   K. Zhou,   Q.  Liu,  J.  Kurths,  Solving Fokker-Planck equation using deep learning. Chaos. 30(1)(2020) 013133.

\bibitem{Zhang2022} H. Zhang,  Y. Xu,   Q.  Liu,  X. Wang,   Y. Li,  Solving Fokker-Planck equations using deep kd-tree with a small amount of data. Nonlinear Dynamics. (2022) https://doi.org/10.1007/s11071-022-07361-2.


\bibitem{li2021data}
Y. Li,   J.  Duan,   A data-driven approach for discovering stochastic dynamical systems with non-Gaussian L\'{e}vy noise. Physica D.  417(2021) 132830.

\bibitem{Lu}Y.Lu,  Y.Li, J. Duan, Extracting stochastic governing laws by non-local Kramers¨CMoyal formulaePhil. Trans. R. Soc. A.380, (2022)20210195.






\bibitem{chen2021solving}
 X. Chen,   L. Yang,   J.  Duan,   G. E. Karniadakis,  Solving Inverse Stochastic Problems from Discrete Particle Observations Using the Fokker-Planck Equation and Physics-Informed Neural Networks. SIAM Journal on Scientific Computing.  43(3)(2021) B811-30.

\bibitem{Yang_I}
Y. Yang, L. Nurbekyan, E. Negrini, R. Martin, M. Pasha. Optimal transport for parameter identification of chaotic dynamics via invariant measures. arXiv preprint arXiv:2104.15138 (2021).



\bibitem{Duan2015} J.  Duan,   \emph{An Introduction to Stochastic Dynamics}, New York(2015) Cambridge University Press.

\bibitem{felix_SDE}
F. Dietrich, A. Makeev, G. Kevrekidis, N. Evangelou, T. Bertalan, S. Reich, I. Kevrekidis. Learning effective stochastic differential equations from microscopic simulations: combining stochastic numerics and deep learning. arXiv preprint arXiv:2106.09004 (2021).


\bibitem{chen_learnMTP}
X. Chen, J. Duan, J. Hu, D. Li. Data-driven method to learn the most probable transition pathway and stochastic differential equation. Physica D: Nonlinear Phenomena 443 (2023) 133559.


\bibitem{Cha}S. Cha,   Comprehensive Survey on Distance/Similarity Measures between Probability Density Functions. \emph{International Journal of Mathematical Models and Methods in Applied Sciences}. 1(4)(2007) 300-307.

\bibitem{Beran}R. Beran,   Minimum hellinger distance estimates for parametric models. Annals of Statistics. 5(3)(1977)445-463.

\bibitem{Klebaner}
  F. C. Klebaner,  \emph{Introduction to  Stochastic  Calculus with Applications}.
 Imperial College Press, London, (2005)2nd edition.



\bibitem{Raissi2019}M. Raissi, P. Perdikaris,   G. E.  Karniadakis,  Physics-informed neural networks: A deep learning framework for solving forward and inverse problems involving nonlinear partial differential equations. Journal of Computational Physics. 378(2019) 686-707.


\bibitem{chen2020}
X. Chen,  J.  Duan, G. E.  Karniadakis,   Learning and meta-learning of stochastic advection-diffusion-reaction systems from sparse measurements. European Journal of Applied Mathematics.  32(3) (2020)397-420.


\bibitem{wang1}
H. Wang,   X. Cheng,  J. Duan,   J. Kurths,  X.  Li,  Likelihood for transcriptions in a genetic regulatory system under asymmetric stable L\'{e}vy noise. Chaos.  28 (2018) 013121.

\bibitem{wang2} X. Cheng,  H. Wang,  X. Wang, J. Duan,  X.  Li, Most probable transition pathways and maximal likely trajectories in a genetic regulatory system.  Physica A. 531 (2019) 121779.



\bibitem{parallel-PINN}
K. Shukla, D. Ameya, G. Karniadakis. Parallel physics-informed neural networks via domain decomposition. Journal of Computational Physics 447 (2021): 110683.

\bibitem{auto}
A. G. Baydin, B. A. Pearlmutter, A. A. Radul, J. M. Siskind. Automatic differentiation in machine learning: a survey. Journal of Marchine Learning Research, 18 (2018): 1-43.

\bibitem{chen2019most}
X. Chen, F. Wu, J. Duan,J. Kurths, X. Li. Most probable dynamics of a genetic regulatory network under stable L{\'e}vy noise. Applied Mathematics and Computation, 348 (2019): 425-436.

\bibitem{Active_learning}
X. Yang, Y. Liu, C. Mi, X. Wang. Active learning Kriging model combining with kernel-density-estimation-based importance sampling method for the estimation of low failure probability. Journal of Mechanical Design, 140 (2018): 051402.
\end{thebibliography}

\end{document}